\documentclass[twoside,a4paper,reqno,11pt]{amsart}
\usepackage{amsfonts, amsbsy, amsmath, amsthm, amssymb, latexsym, verbatim, enumerate}
\usepackage{mathrsfs}
\usepackage{bm}
\usepackage[pdftex]{color,graphicx}

\makeatletter
\newcommand{\imod}[1]{\allowbreak\mkern4mu({\operator@font mod}\,\,#1)}
\makeatother

\def\hal{\unskip\nobreak\hfill\penalty50\hskip10pt\hbox{}\nobreak

\hfill\vrule height 5pt width 6pt depth 1pt\par\vskip 2mm}

\usepackage{bbm}
\usepackage[all,cmtip]{xy}

\usepackage{enumerate}
\usepackage[shortlabels]{enumitem}

 \def\Z{\mathbb Z}
  \def\N{\mathbb N}
  \def\F{\mathbb F}

\def\i{\iota}
\def\c{\chi}
  \def\a{\alpha}
  
  \def\g{\gamma}
  \def\G{\Gamma}
\def\d{\delta}
   \def\D{\Delta}
  \def\e{\epsilon}
  \def\la{\lambda}

  \def\l{\langle}
  \def\r{\rangle}
  \def\s{\sigma}

    \def\O{\Omega}
  \def\no{\noindent}
  \def\go{\rightarrow}
  \def\pf{\noindent {\bf Proof $\;$ }}

\newcommand{\Irr}{{\mathrm {Irr}}}

\newcommand{\Hom}{{\mathrm {Hom}}}

\newcommand{\CB}{{\bf C}}

\newcommand{\NB}{{\bf N}}

\newcommand{\KK}{{\mathbb K}}

\newcommand{\GC}{{\mathcal G}}

\newcommand{\LC}{{\mathcal L}}

\newcommand{\EC}{{\mathcal E}}

\renewcommand{\mod}{\bmod \,}

\numberwithin{equation}{section}

  \def\hal{\unskip\nobreak\hfil\penalty50\hskip10pt\hbox{}\nobreak
  \hfill\vrule height 5pt width 6pt depth 1pt\par\vskip 2mm}

\begin{document}


\author[M. W. Liebeck]{Martin W. Liebeck}
\address{M.W. Liebeck, Department of Mathematics,
    Imperial College, London SW7 2BZ, UK}
\email{m.liebeck@imperial.ac.uk}

\author[A. Shalev]{Aner Shalev}
\address{A. Shalev, Institute of Mathematics, Hebrew University, Jerusalem 91904, Israel}
\email{shalev@math.huji.ac.il}

\author[P. H. Tiep]{Pham Huu Tiep}
\address{P.H. Tiep, Department of Mathematics, Rutgers University, Piscataway, NJ 08854, USA}
\email{tiep@math.rutgers.edu}

\title[Character ratios, representation varieties and generation]{Character ratios, representation varieties and random generation of finite groups of Lie type}

  \maketitle

\begin{abstract}
We use character theory of finite groups of Lie type to establish new results on representation varieties of Fuchsian groups, and also on probabilistic generation of groups of Lie type.
\end{abstract}

\footnotetext{The first and second authors acknowledge the support of EPSRC grant EP/H018891/1.}
\footnotetext{The second author acknowledges the support of ERC advanced grant
247034, ISF grants 1117/13 and 686/17, BSF grant 2016072,  and the Vinik chair of mathematics which he holds.}
\footnotetext{The third author was partially supported by the NSF grant DMS-1839351 and DMS-1840702, the Simons Foundation
Fellowship 305247, the EPSRC, the Mathematisches Forschungsinstitut Oberwolfach, and the Joshua Barlaz Chair in Mathematics.
Parts of the paper were written while the third author visited the Departments of Mathematics of Imperial College London, and
University of Chicago. It is a pleasure to thank Imperial College, Prof. Ngo Bao Chau, and the University of Chicago
for generous hospitality and stimulating environment.}
\footnote{The authors are grateful to the referee for insightful comments on the paper.}

  \newtheorem{thm}{Theorem}[section]
  \newtheorem{prop}[thm]{Proposition}
  \newtheorem{lem}[thm]{Lemma}
  \newtheorem{cor}[thm]{Corollary}
   \newtheorem{rem}[thm]{Remark}
   \newtheorem{exa}[thm]{Example}

\tableofcontents


\section{Introduction}

The purpose of this paper is to apply some new character-theoretic results on finite groups of Lie type proved in \cite{blst1} to the study of representation varieties of Fuchsian groups, and also to probabilistic generation of finite simple groups. Our results represent substantial advances on those obtained on these topics in \cite{fuchs2}
and \cite{LL}.

Recall that a Fuchsian group is a finitely generated non-elementary discrete
group of isometries of the hyperbolic plane. These groups play an important role in geometry,
analysis and algebra; they have been studied extensively since the days of Fricke and Klein,
who found presentations for them, which will be described below.
Important examples of Fuchsian groups include free groups, the modular group, surface groups,
the Hurwitz group and hyperbolic triangle groups in general.

In the past two decades there has been considerable renewed interest in Fuchsian groups and
their homomorphisms to finite simple groups. This was partly motivated by Higman's conjecture
that every Fuchsian group surjects to all large enough alternating groups, which was proved by
Everitt \cite{Ev} (see also Conder \cite{C}) in the oriented case. A probabilistic proof of
Higman's conjecture which includes the non-oriented Fuchsian groups appeared in \cite{fuchs1}.

A second motivation to study homomorphisms from Fuchsian groups to finite simple groups
is the attempt to generalize various results on random generation of finite simple groups,
for example by two elements \cite{lish95}, by elements of orders $2$ and $3$ \cite{lish96}, and so on.
These results show that a random homomorphism from the free group $F_2$, or from the modular group
$PSL_2(\Z)$, is an epimorphism (with some exceptions in the latter case). Can one obtain similar results
for general Fuchsian groups?

The answer is positive if we study homomorphisms to alternating groups, but less is known
for homomorphisms from Fuchsian groups to simple groups of Lie type. In \cite{fuchs2} it is
shown that random homomorphisms from oriented Fuchsian groups of genus at least 2 (and from
non-oriented Fuchsian groups of genus at least 3) are onto. The remaining case of Fuchsian groups
of very small genus has resisted all attacks so far.

Note that some Fuchsian groups do not surject to all large enough finite simple groups (for example
the simple groups may not have elements of prescribed orders).
Important results on Fuchsian generation of finite simple groups, focusing on hyperbolic triangle
groups (which have genus zero) were obtained by Larsen, Lubotzky and Marion \cite{LLM1}, \cite{LLM2}
using deformation theory. However, they do not cover random Fuchsian (or triangle) generation.

In \cite[p.323]{fuchs2} it was conjectured that a random homomorphism from any Fuchsian group to a finite simple classical group of large enough rank is onto. This paper establishes various cases of this conjecture.

We are also interested in representation varieties of Fuchsian groups and in studying their dimensions.
By counting homomorphisms of Fuchsian groups to finite classical groups and using Lang-Weil estimates we draw
conclusions concerning representation varieties.

Again, satisfactory results were obtained in \cite{fuchs2} for Fuchsian groups of genus at least 2
(at least 3 in the non-oriented case). Asymptotic results over fields of characteristic zero were
subsequently obtained using deformation theory \cite{LL}, but the general case remained open.
In this paper we manage to solve various instances of this challenging problem.
Strong character bounds for finite groups of Lie type, recently obtained in \cite{blst1},
play a major role in our proofs.

We now describe our results in detail.
Let $\G$ be a co-compact Fuchsian group of genus $g$ having $d$ elliptic generators of orders $m_1,\ldots ,m_d$ (all at least 2). Thus if $\G$ is  orientation-preserving, it has a presentation of the following form:
\[
\l a_1,b_1,\ldots ,a_g,b_g,\; x_1,\ldots ,x_d \;|\;x_1^{m_1} = \cdots = x_d^{m_d} = 1, \;
 x_1\cdots x_d \, \prod_{i=1}^g[a_i,b_i] = 1 \r,
\]
and if $\G$ is non-orientation-preserving it has  a presentation
\[
\l a_1,\ldots ,a_g,\; x_1,\ldots ,x_d\;|\;    x_1^{m_1} = \cdots = x_d^{m_d} = 1, \;
 x_1\cdots x_d \,a_1^2\cdots  a_g^2 = 1\r
\]
(see \cite{beardon}). The {\it measure} of $\G$ is
\begin{equation}\label{meas}
\mu(\G) = vg-2+\sum_{i=1}^d \left(1-\frac{1}{m_i}\right) > 0,
\end{equation}
where $v = 2$ if $\G$ is orientation-preserving and $v=1$ if not.
We assume throughout that $\G$ is not virtually abelian, which implies that $vg+d\ge 3$.

Write $\mu = \mu(\G)$. We will prove various results showing that (under some conditions),
the number of homomorphisms from a Fuchsian group $\G$ to a classical group $G = G_n(q)$
is roughly $|G|^{\mu + 1}$, up to some multiplicative error terms.
These error terms are relatively small when the dimension $n$ is large and when $q$ is
much larger than $n$.
We will then conclude that the dimension of the respective
representation variety $\Hom(\Gamma, G_n(\KK))$ for algebraically closed fields
$\KK$ is roughly $(\mu + 1)\dim G_n(\KK)$, up to a small explicit additive error term.

Define
\[
N_1(\G) = {\rm max} \left(\frac{2+\sum \frac{1}{m_i}}{\mu},m_1,\ldots,m_d\right)+1.
\]
If $\mu>2$, define also
\[
N_2(\G) = {\rm max} \left(N_1(\G),\, \frac{d+16}{4(\mu -2)}+1\right).
\]

\begin{thm}\label{fuchsian}
There exist functions $f,g:\N \to \N$ such that the following statement holds. Suppose that $n\ge N_1(\G)$, $q > f(n)$,  and $q \equiv 1 (\bmod\ m_i)$ for all
$i$. Then
\begin{itemize}
\item[{\rm (i)}] $|{\rm Hom}(\G,SL_n(q))| > |SL_n(q)|^{\mu+1}q^{-\sum m_i}$.
\item[{\rm (ii)}]  If $\mu > 2$ and also $n\ge N_2(\G)$, then
$$|{\rm Hom}(\G,GL_n(q))| < g(n)\,q\, |GL_n(q)|^{\mu+1}.$$
\end{itemize}
\end{thm}

This follows from \cite{fuchs2} if the genus $g$ of $\G$ is at least 2 (at
least 3 if $\G$ is non-oriented), but is new for groups of genus 0 or 1,
which are the most challenging cases.

For $\KK$ an algebraically closed field, the representation variety of $\G$ in dimension $n$ over $\KK$ is defined to be
$R_{n,\KK}(\G): = {\rm Hom}(\G,GL_n(\KK))$. Using Theorem \ref{fuchsian} we can get the following quite precise estimates for the dimension of this variety.

\begin{cor}\label{dimhomgl}
Suppose $n \ge N_1(\G)$, and let $\KK$ be an algebraically closed field of characteristic not dividing $m_1\cdots m_d$.
\begin{itemize}
\item[{\rm (i)}] Then $\dim R_{n,\KK}(\G) \ge (n^2-1)(1+\mu)-\sum_{i=1}^d m_i$.
\item[{\rm (ii)}]  If also $\mu > 2$ and $n\ge N_2(\G)$, then
\[
\dim R_{n,\KK}(\G) = n^2(1+\mu) - c,
\]
where $-1 \le c \le \mu + 1 +\sum_{i=1}^d m_i$.
\item[{\rm (iii)}] The upper bound for $|{\rm Hom}(\G,GL_n(q))|$ given in Theorem {\rm \ref{fuchsian}(ii)} holds for all $q$ coprime to $m_1\cdots m_d$. 
\end{itemize}
\end{cor}

For genus $g\ge 2$ ($g\ge 3$ if $\G$ is non-oriented), the exact value of $\dim R_{n,\KK}(\G)$ is given by \cite[1.9]{fuchs2}, but the above result covers all possible $g$.

Note that \cite{LL} contains asymptotic results for the dimensions of representation varieties of all Fuchsian groups in characteristic 0.


\vspace{4mm}
We also prove analogous results for other classical groups. To state these, let $\G$ be a Fuchsian group as above, suppose $m_i$ is odd for all $i$, and define
\[
t(\G) = \sum_{i=1}^d \left(\frac{1}{m_i-1}-\frac{1}{m_i}\right).
\]
If $\mu > t(\G)$, define
\[
N_3(\G) = {\rm max} \left(\frac{1+\sum \frac{2}{m_i-1}}{\mu- t(\G)},\frac{m_1^3}{3},\ldots,\frac{m_d^3}{3}\right),
\]
and if $\mu > 2$, define
\[
N_4(\G) =  {\rm max} \left(\frac{d+8}{4(\mu-2)},m_1^2,\ldots ,m_d^2\right).
\]

\begin{thm}\label{fuchsclass}
Suppose that $m_i$ is odd for all $i$, and let $G_n(q)$ be one of the classical groups $Sp_{2n}(q)$, $\O_{2n+1}(q)$, $\O_{2n+2}^\pm(q)$, where $q \equiv 1 (\mod 2m_i)$ for all $i$. Then provided $q>f(n)$, the following hold (for suitable functions $f,g$).
\begin{itemize}
\item[{\rm (i)}] If $\mu > t(\G)$ and $n>N_3(\G)$, then
\[
|{\rm Hom}(\G,G_n(q))| > |G_n(q)|^{\mu+1}q^{-\sum \frac{m_i^2}{2}}.
\]
\item[{\rm (ii)}]  If $\mu > 2$ and $n > N_4(\G)$, then
\[
|{\rm Hom}(\G,G_n(q))| < g(n)\, q^d\,|G_n(q)|^{\mu+1}.
\]
\end{itemize}
\end{thm}

Note that the assumption that the $m_i$ are odd is needed so that in the proof, elements of order $m_i$ can be chosen to satisfy the splitness condition required for the application of the character-theoretic bound given in Theorem \ref{charbd}. A similar comment applies to the assumptions in Theorems \ref{excep1} -- \ref{exceptri} below.

As before, Theorem \ref{fuchsclass}  leads to corresponding estimates for the dimensions of the representation varieties ${\rm Hom}(\G,G_n(\KK))$ for an algebraically closed field $\KK$:

\begin{cor}\label{dimhomclass}
Suppose that all the $m_i$ are odd, $\mu > {\rm max}\left(2,\,t(\G)\right)$ and also $n> {\rm max}\left(N_3(\G),N_4(\G)\right)$, and
let $G_n(q)$ be as in Theorem $\ref{fuchsclass}$. 
\begin{itemize}
\item[{\rm (i)}] For an algebraically closed field $\KK$ of characteristic not dividing \linebreak $2m_1\cdots m_d$,
\[
 \dim {\rm Hom}(\G, G_n(\KK)) = (\mu+1) \dim G_n(\KK)-c,
\]
where $-d \le c \le \frac{1}{2}\sum m_i^2$.
\item[{\rm (ii)}] The upper bound for $|{\rm Hom}(\G,G_n(q))|$ given in Theorem {\rm \ref{fuchsclass}(ii)} holds for all $q$ coprime to $2m_1\cdots m_d$. 
\end{itemize}
\end{cor}

We can use the above results to deduce consequences on random Fuchsian generation of classical groups. In the next result, for a finite group $X$ denote by $P_{\G}(X)$ the probability that a randomly chosen homomorphism in ${\rm Hom}(\G,X)$ is an epimorphism.

For $SL_n(q)$ we prove the following.

\begin{thm}\label{fuchsgensl}
Let $\G$ be a Fuchsian group as above, and assume
\[
\mu > \nu:= {\rm max}\left(2,\,1+\sum\frac{1}{m_i}\right).
\]
Define
\[
Q = \bigcup_{primes\;p} \{q : q=p^a \equiv 1 (\mod m_i)\;\forall i\}.
\]
Then for $n \ge \nu N_2(\G)+2\sum m_i$, we have
\[
{\rm lim}_{q\rightarrow \infty,\,q \in Q} P_\G(SL_n(q)) = 1.
\]
\end{thm}

This follows from \cite[Theorem 1.6]{fuchs2} for $g\ge 2$  ($g\ge 3$ if $\G$ is non-oriented) without any assumptions on $\mu$ or $n$. But it is new for $g\le 1$ -- indeed, it partially proves the conjecture in \cite[p.323]{fuchs2}.

Taking $g=0$, we can deduce the following generation result from
Theorem \ref{fuchsgensl}. Recall that a group $G$ is said to be $(m_1,\ldots ,m_d)$-generated
if it is generated by elements $g_1, \ldots , g_d$ satisfying
$g_1^{m_1} = \ldots = g_d^{m_d} = g_1 \cdots g_d = 1$.

\begin{cor}\label{slngen}
Suppose $d\ge 5$ and $m_1,\ldots ,m_d$ are integers, all at least $2$, such that
\[
d > {\rm max} \left(3+2\sum_{i=1}^d \frac{1}{m_i},\; 4+\sum_{i=1}^d \frac{1}{m_i} \right).
\]
Define $\nu, N_2(\G)$ as above. Then for $n \ge  \nu N_2(\G)+\sum m_i$, and for all sufficiently large $q$ such that $q \equiv 1 (\mod m_i)$ for all $i$, the group $SL_n(q)$ is $(m_1,\ldots ,m_d)$-generated.
\end{cor}


We can prove results of a similar flavour for the other classical groups. Let $\G$ be a Fuchsian group as above, and assume the following:
\begin{itemize}
\item[(a)] all the $m_i$ are odd, and
\item[(b)] $\mu > {\rm max} \left(2,\, t(\G),\,1+\sum \frac{1}{m_i}\right)$.
\end{itemize}
Write $\s_1 = \sum \frac{1}{m_i}$, $\s_2 = \sum m_i$, $\s_3 = \frac{1}{2}\sum m_i^2$, and define
\[
N_5(\G) = {\rm max} \left( N_3(\G),\,2N_2(\G),\, \s_1 N_2(\G)+\s_3 + 2,\, 2N_4(\G),\,
(1+2\s_1)N_4(\G)+\s_2   \right).
\]

\begin{thm}\label{fuchsgenclass}
Let $\G$ be a Fuchsian group, and assume conditions $(a)$ and $(b)$ above hold.
Let $G_n(q)$ be $Sp_{2n}(q)$, $\O_{2n+1}(q)$, or $\O_{2n+2}^\pm(q)$ with $n> N_5(\G)$, and
define
\[
Q = \bigcup_{primes\;p} \{q : q=p^a \equiv 1 (\mod 2m_i)\;\forall i\}.
\]
Then
\[
{\rm lim}_{q\rightarrow \infty,\,q \in Q} P_\G(G_n(q)) = 1.
\]
\end{thm}

A deduction similar to Corollary \ref{slngen} can be made concerning the $(m_1,\ldots ,m_d)$-generation of $G_n(q)$; we leave it to the interested reader to
work this out.

\vspace{4mm}
Finally, we turn to versions of the above results for exceptional groups, provided the elliptic orders $m_i$ are not too small. In this case we are able to provide a precise formula for the dimension of the respective
representation variety.
For an algebraic group $\GC$ and a positive integer $m$, let $J_m(\GC)$ be the variety $\{x \in \GC : x^m=1\}$. The dimensions of these varieties are given by \cite{law}.

\begin{thm}\label{excep1}
Let $\G$ be a Fuchsian group as above, and
let $\GC$ be a simple adjoint algebraic group of exceptional type over an algebraically closed
field $\KK$ of good characteristic. Suppose that $m_1\cdots m_d$ is coprime to $30$ and also to ${\rm char }(\KK)$.
 Then
\[
\dim {\rm Hom}(\G,\GC) = (vg-1)\dim \GC + \sum_{i=1}^d \dim J_{m_i}(\GC).
\]
\end{thm}

See Theorem \ref{class1} for a variant of this result for classical groups of small rank.

We also obtain results on random Fuchsian generation of exceptional groups of Lie type.

\begin{thm}\label{excep2} Let $\G$ be a Fuchsian group as above, and
let $G(q)$ denote a simple group of exceptional Lie type over $\F_q$. Suppose that $m_1\cdots m_d$ is coprime to $30$, and
define
\[
Q = \bigcup_{good\;primes\;p}\{q : q=p^a \equiv 1 (\mod m_i) \;\forall i\}.
\]
Then
\[
{\rm lim}_{q\rightarrow \infty,\,q \in Q} P_\G(G(q)) = 1.
\]
In particular, for all sufficiently large such $q$, $G(q)$ is a quotient of $\G$.
\end{thm}

Applying this in the important and challenging case where $\G$ is a triangle group, we obtain the following
random triangle generation theorem for exceptional groups.

\begin{thm}\label{exceptri}  Let $m_1,m_2,m_3$ be positive integers that are coprime to $30$ and let $p$ be a
good prime for the exceptional algebraic group $\GC$ that does not divide $m_1m_2m_3$. Let $q = p^a \equiv 1 (\mod m_i)$
for all $i$. Then, as $q \to \infty$, the group $G(q) = \GC^{F_q}$ is randomly $(m_1,m_2,m_3)$-generated.
In particular, for all sufficiently large such $q$, $G(q)$ is $(m_1,m_2,m_3)$-generated.
\end{thm}

Our results above on random Fuchsian generation and random triangle generation are novel.

A variety of important triangle generation results for finite simple groups $G(q)$ of Lie type were obtained in
\cite{LLM1, LLM2} using deformation theory. These results show that (under some assumptions on $m_1, m_2, m_3$
and the Lie type $G$), there are integers $p_0, e > 0$ such that

(i) for all primes $p > p_0$ and for all integers $\ell \ge 1$, the group $G(p^{e \ell})$ is $(m_1,m_2,m_3)$-generated, and

(ii) the set of primes $p$ for which $G(p^{\ell})$ is $(m_1,m_2,m_3)$-generated for all $\ell \ge 1$ has positive density.

While the triangle generation results of \cite{LLM1, LLM2} are more general than ours, they are
non-effective. Consequently, the numbers $p_0$ and $e$ cannot be given explicitly. 

The generation result in Theorem \ref{exceptri} above has several advantages. First, it implies condition (i) with explicit $p_0 := \max (m_1, m_2, m_3)$ (with implicit $e$).
Secondly, setting $m = {\mathrm {lcm}}(m_1, m_2, m_3)$ and $a = \phi(m)$, the Euler function of $m$, it implies
that the group $G(p^{a \ell})$ is $(m_1,m_2,m_3)$-generated for all sufficiently large $\ell$ (a condition stronger
than condition (i) with implicit $e$).
Thirdly, it shows that the set of primes $p$ satisfying condition (ii) includes all sufficiently large primes in the arithmetic progression $\{ im + 1: i \ge 1 \}$ and so has density at least $1/a$.

\section{Character bounds} \label{charbds}

In this section we summarize the results in \cite{blst1} on character ratios for groups of Lie type that we will need in our proofs.

Let $\KK$ be an algebraically closed field of characteristic $p>0$,
$\GC$ a connected reductive algebraic group over $\KK$, $F:\GC \to \GC$ a Frobenius endomorphism,
and $G = \GC^F$.
For a subgroup $L$ of $\GC$ write $L_{{\rm {unip}}}$ for the set of non-identity unipotent elements
of $L$. For a fixed $F$, a Levi subgroup $\LC$ of $\GC$ will be called {\it split}, if it is an $F$-stable
Levi subgroup of an $F$-stable proper parabolic subgroup of $\GC$. We define
\[
\a(\LC) = {\rm max}_{u \in \LC_{{\rm unip}}} \frac{\dim u^\LC}{\dim u^\GC},
\]
if $\LC$ is not a torus, and $\a(\LC) = 0$ otherwise.

\begin{thm}\label{charbd} {\rm (\cite[Theorem 1.1]{blst1})} There exists a function $f:\N \to \N$ such that the
following statement holds.
With the above notation, assume $[\GC,\GC]$ is simple of rank $r$ and $p$ is a good prime for $\GC$.
Let $g \in  G$ be any element such that $\CB_G(g) \le L:= \LC^F$, where $\LC$ is a split Levi subgroup of $\GC$.
Then for any irreducible character $\c$ of $G$,
\[
|\c(g)| \le f(r)\,\c(1)^{\a(\LC)}.
\]
\end{thm}

We note that the function $f$ above is explicit, see \cite[Remark 1.2]{blst1}.

In order to apply Theorem \ref{charbd}, one needs information about the values $\a(\LC)$ for Levi subgroups $\LC$.
For classical groups the following bound is proved in \cite{blst1}.

\begin{thm} \label{ratio} {\rm (\cite[Theorem 1.4]{blst1})} If $\GC$ is a classical algebraic group in good characteristic, and $\LC$ is a Levi subgroup of $\GC$, then $\a(\LC) \le \frac{1}{2}\left(1+\frac{\dim \LC}{\dim \GC}\right)$.
\end{thm}

Some more precise bounds for the values of $\a(\LC)$ for specific types of Levi subgroups of classical groups can be found in the proofs of Propositions \ref{fomin}--\ref{fominclass} below. For exceptional groups we have the following.

\begin{thm} \label{alphaexcep} {\rm (\cite[Theorem 1.4]{blst1})} If $\GC$ is an exceptional algebraic group in  good characteristic, the values of $\a(\LC)$ for Levi subgroups $\LC$ are as in Table $\ref{extab}$.
\end{thm}

In Table \ref{extab}, for $\GC = F_4$ or $G_2$ the symbols $\tilde A_1, \tilde A_2$ refer to Levi subsystems consisting of short roots. For $\GC = E_7$, there are two Levi subgroups $A_5$ and $A_5'$: using the notation for the fundamental roots $\a_i\,(1\le i\le 7)$ as in \cite{bour}, these are the Levi subgroups with fundamental roots $\{\a_i : i=1,3,4,5,6\}$ and $\{\a_i : i=2,4,5,6,7\}$ respectively. The notation  $\triangleright A_4$, for instance, means that
$\LC'=[\LC,\LC]$ has a simple factor of type $A_4$.

\begin{table}[h!]
\caption{$\a$-values for exceptional groups} \label{extab}
\vspace{-7mm}
\[
\begin{array}{r|cccccccccc}
\hline
&&&&&&&&&& \\
\vspace{-8mm}\\
\GC=E_8,\,\LC '= & E_7 & D_7 & \LC'\triangleright E_6 & D_6 & A_7 & \triangleright D_5 & \triangleright A_6 & \triangleright A_5 &
\triangleright D_4 & \hbox{rest} \\
\a(\LC)= & \frac{17}{29} & \frac{9}{23} & \frac{11}{29} & \frac{9}{29} & \frac{15}{56} &  \frac{7}{29} & \frac{5}{23}
&  \frac{4}{23} &  \frac{5}{29} &  \le \frac{1}{6}  \\
\vspace{-3.5mm}\\
\hline
&&&&&&&&&& \\
\vspace{-8mm}\\
\GC=E_7,\,\LC '= & E_6 & D_6 & \LC'\triangleright D_5 & A_6 & A_5 & \triangleright A_5' & \triangleright D_4 & \triangleright A_4 &
\triangleright A_3 & \hbox{rest} \\
\a(\LC)= & \frac{11}{17} & \frac{5}{9} & \frac{7}{17} & \frac{5}{13} & \frac{4}{13} &  \frac{1}{3} & \frac{5}{17}
& \le \frac{1}{4} &  \le \frac{1}{5} &  \le \frac{1}{6}  \\
\vspace{-3.5mm}\\
\hline
&&&&&&&&&& \\
\vspace{-8mm}\\
\GC=E_6,\,\LC '= & D_5 & A_5 & D_4 & \LC'\triangleright A_4 & \triangleright A_3 & \triangleright A_2 & A_1^k &&&\\
\a(\LC)= & \frac{7}{11} & \frac{1}{2} & \frac{5}{11} & \frac{3}{8} & \frac{3}{11} & \le \frac{7}{27} & \le \frac{3}{20}&&& \\
\vspace{-3mm}\\
\hline
&&&&&&&&&& \\
\vspace{-8mm}\\
\GC=F_4,\,\LC ' =& B_3 & C_3 & A_2\tilde A_1,A_2 & \tilde A_2A_1 & \tilde A_2 & A_1\tilde A_1 & A_1 & \tilde A_1 &&\\
\a(\LC)= & \frac{1}{2} & \frac{7}{15} & \frac{1}{4} & \frac{2}{9} & \frac{1}{5} & \frac{1}{7} & \frac{1}{8} & \frac{1}{11}&& \\
\vspace{-3.5mm}\\
\hline
&&&&&&&&&& \\
\vspace{-8mm}\\
\GC=G_2,\,\LC '= & A_1 & \tilde A_1&&&&&&&& \\
\a(\LC)= & \frac{1}{3} & \frac{1}{4}&&&&&&&& \\
\hline
\end{array}
\]
\end{table}

For calculating $\a$-values in $SL_n(\KK)$, the following elementary result is sometimes useful.

\begin{lem}\label{sscomp} {\rm (\cite[Cor. 4.2]{blst1})}
If $\LC$ is a Levi subgroup of $\GC = GL_n(\KK)$ or $SL_n(\KK)$, then, denoting by $\LC_{ss}$ the set of
non-central semisimple elements in $\LC$,
\[
\a(\LC) \leq \hbox{{\rm {max}}}_{s \in \LC_{ss}} \frac{\dim s^{\LC}}{\dim s^{\GC}}.
\]
\end{lem}

We will also require the next result, taken from \cite{chardeg}.
Suppose $\GC$ is a  simple algebraic group of rank $r$ in good characteristic,
and $G = G(q)= \GC^F$ is a finite quasisimple group of Lie type over $\F_q$.
For a real number $s$, define
\begin{equation}\label{zetadef}
\zeta^G(s) = \sum_{\c \in \Irr(G)} \c(1)^{-s}.
\end{equation}

\begin{thm}\label{zetabd} {\rm (\cite[Theorem 1.1]{chardeg})}
Let $h$ be the Coxeter number of $\GC$.
If $s>\frac{2}{h}$, then $\zeta^G(s) \rightarrow 1$ as $q\rightarrow \infty$.
\end{thm}

\section{Fuchsian groups and representation varieties
}

In this section we prove Theorems
\ref{fuchsian}, \ref{fuchsclass} and \ref{excep1}, and also deduce Corollaries  \ref{dimhomgl} and \ref{dimhomclass}.
For the proofs we need some sharp character bounds on specific elements of classical groups, which are given in the first subsection.

\subsection{Some sharp character bounds}

The first result gives a sharp bound on the character values of particular elements of $GL_n(q)$.

\begin{prop}\label{fomin}
Let $n\ge m\ge 2$ be integers, and write $n=km+s$ with $0\le s<m$. Let $q$ be a prime power such that $q \equiv 1(\mod m)$.
Then the group $GL_n(q)$ contains an  element $x$ with the following properties:
\begin{itemize}
\item[{\rm (i)}] $x$ has order $m$;
\item[{\rm (ii)}] $\CB_{GL_n(q)}(x) = GL_{k+1}(q)^s \times GL_k(q)^{m-s}$;
\item[{\rm (iii)}] for all $\chi \in \Irr(GL_n(q))$, $|\chi(x)| < f(n)\,\chi(1)^{\frac{1}{m}}$.
\end{itemize}
\end{prop}

An adjustment is needed if we require the element $x$ to lie in $SL_n(q)$:

\begin{prop}\label{fominsl}
Let $n,m,k,s$ and $q$ be as in Proposition $\ref{fomin}$.
Then the group $SL_n(q)$ contains an  element $x$ with the following properties:
\begin{itemize}
\item[{\rm (i)}] $x$ has order $m$;
\item[{\rm (ii)}] we have
\[
\CB_{GL_n(q)}(x) = \left\{\begin{array}{ll}
GL_{k+1}(q) \times GL_{k-1}(q)\times GL_k(q)^{m-2}, & \hbox{if }n=km \hbox{ with } \\
                                                                     & m \hbox{ even}, k\hbox{ odd} \\
GL_{k+1}(q)^s \times GL_k(q)^{m-s}, \hbox{ otherwise} &
\end{array}
\right.
\]
\item[{\rm (iii)}] for all $\chi \in \Irr(SL_n(q))$,
\[
|\chi(x)| <
\left\{\begin{array}{l}
f(n)\,\chi(1)^{\frac{1}{m}\left(1+\frac{1}{n-1}\right)},  \hbox{ if }n=km \hbox{ with }m \hbox{ even}, k\hbox{ odd} \\
f(n)\,\chi(1)^{\frac{1}{m}}, \hbox{ otherwise}
\end{array}
\right.
\]
\end{itemize}
\end{prop}


We also prove analogous results for other classical groups:

\begin{prop}\label{fominclass}
Let $n,m$ be integers with $m$ odd, $m\ge 3$ and $n>\frac{m^3}{3}$. Suppose $q \equiv 1 (\mod 2m)$, and let $G=G_n(q)$ be one of the classical groups $Sp_{2n}(q)$, $\O_{2n+1}(q)$, $\O_{2n+2}^\pm (q)$. Then $G_n(q)$ contains an element $x$ with the following properties:
\begin{itemize}
 \item[{\rm (i)}] $x$ has order $m$;
 \item[{\rm (ii)}] $|x^G|>c|G|^{1-\frac{1}{m}}q^{-\frac{m^2}{2}}$, where $c$ is a positive absolute constant;
 \item[{\rm (iii)}] for all $\chi \in \Irr(G_n(q))$,
\[
|\chi(x)| < f(n)\,\chi(1)^{\frac{1}{m-1}\left(1+\frac{2}{n}\right)}.
\]
\end{itemize}
\end{prop}

\vspace{6mm}
\no {\bf Proof of Proposition \ref{fomin}}

\vspace{2mm}
Let $n,m$ be integers with $n\ge m\ge 2$, and
let $p$ be a prime and $q= p^a \equiv 1 (\mod m)$. Define $\GC = GL_n(\KK)$ where $\KK = \bar \F_p$, and let $F: \GC \rightarrow \GC$ be the Frobenius endomorphism sending $(a_{ij}) \rightarrow (a_{ij}^q)$, so that $G:= \GC^F = GL_n(q)$.
Let $\la_1,\ldots ,\la_m$ be the $m^{th}$ roots of 1 in $\F_q^*$. Write $n=km+s$ with $0\le s<m$, and define
\begin{equation}\label{xdefn}
x = {\rm diag}\left(\la_1 I_{k+1},\ldots ,\la_s I_{k+1},\la_{s+1} I_{k},\ldots ,\la_m I_{k}\right).
\end{equation}
Then $x$ has order $m$ and  $\CB_G(x) = GL_{k+1}(q)^s \times GL_{k}(q)^{m-s} = \LC^F$, where $\LC = GL_{k+1}(\KK)^s \times GL_{k}(\KK)^{m-s}$ is an $F$-stable split Levi subgroup in $\GC$.

We claim that
\begin{equation}\label{1kineq}
\a(\LC) \le \frac{1}{m}.
\end{equation}
Given this,  the conclusion of Proposition \ref{fomin} follows immediately from Theorem \ref{charbd}.

To prove (\ref{1kineq}) we apply Lemma \ref{sscomp}. We have $\LC = \prod_{i=1}^m GL(V_i)$, where
$\dim V_i = k+1$ for $1\le i \le s$ and $\dim V_i = k$ for $s+1 \le i \le m$. Let $y \in \LC$ be a semisimple element, taken to be a diagonal element in each factor $GL(V_i)$. Let $\mu_1,\ldots ,\mu_r$ be the distinct eigenvalues of $y$ on $V$, and  for each $i$ write
\begin{equation}\label{yab}
y^{V_i} = {\rm diag}\left(\mu_1 I_{a_{i1}},\ldots ,\mu_r I_{a_{ir}}\right),
\end{equation}
so that $\sum_j a_{ij}$ is equal to $k+1$ for $1\le i\le s$, and is equal to $k$ for $s+1\le i\le m$. Then
\[
\begin{array}{l}
\dim y^\LC = s(k+1)^2+(m-s)k^2-\sum_{j=1}^r\sum_{i=1}^m a_{ij}^2, \\
\dim y^\GC = (km+s)^2 -  \sum_{j=1}^r\left(\sum_{i=1}^m a_{ij}\right)^2.
\end{array}
\]
Hence
\[
\begin{array}{lll}
\frac{\dim y^\LC}{\dim y^{\GC}} \le \frac{1}{m} & \Leftrightarrow & m\dim y^\LC \le \dim y^\GC \\
& \Leftrightarrow & ms(k+1)^2+m(m-s)k^2-m\sum_{j=1}^r\sum_{i=1}^m a_{ij}^2 \le \\
&&  (km+s)^2- \sum_{j=1}^r\left(\sum_{i=1}^m a_{ij}\right)^2 \\
& \Leftrightarrow & s(m-s) \le \sum_{j=1}^r \left(m\sum_{i=1}^m a_{ij}^2 - \left(\sum_{i=1}^m a_{ij}\right)^2\right).
\end{array}
\]
The sum in the right hand side of the above inequality is equal to a sum of squares of the form $a_{ij}^2$
and $(a_{ij}-a_{i'j})^2$, and at least
$s(m-s)$ of the latter squares, with $1\le i\le s$ and $s+1 \le i'\le m$, are nonzero. Hence the above inequality holds, and so
$\frac{\dim y^\LC}{\dim y^{\GC}} \le \frac{1}{m}$. By Lemma \ref{sscomp}, this establishes (\ref{1kineq}), completing the proof of Proposition \ref{fomin}. \hal

\vspace{6mm}

\no {\bf Proof of Proposition \ref{fominsl}}

\vspace{2mm}
Let $G = SL_n(q) = \GC^F$, where $\GC = SL_n(\KK)$.
If $n$ is not of the form $km$ with $m$ even, $k$ odd, then we can choose $x$ as in (\ref{xdefn}) to have determinant 1, and the proof goes through exactly as for Proposition \ref{fomin} above.

So assume $n=km$ with $m$ even, $k$ odd. In this case, taking $\la_1=1$, $\la_2 = -1$ we define
\[
x = {\rm diag}\left(I_{k+1},-I_{k-1}, \la_3 I_{k},\ldots ,\la_m I_{k}\right).
\]
Then $x$ has order $m$ and determinant 1, and
$$\CB_G(x) = (GL_{k+1}(q) \times GL_{k-1}(q)\times GL_k(q)^{m-2})\cap G = \LC^F,$$
where $\LC = (GL_{k+1}(\KK) \times GL_{k-1}(\KK)\times GL_k(\KK)^{m-2})\cap \GC$ is an $F$-stable split Levi subgroup of $\GC$. So the proposition will follow from Theorem \ref{charbd} once we prove
\begin{equation}\label{alin}
\a(\LC) \le \frac{1}{m}\left(1+\frac{1}{n-1}\right).
\end{equation}
Let $y \in \LC$ be as above in (\ref{yab}). Then
\[
\begin{array}{l}
\dim y^\LC = mk^2+2-\sum_{j=1}^r\sum_{i=1}^m a_{ij}^2, \\
\dim y^\GC = m^2k^2 -  \sum_{j=1}^r\left(\sum_{i=1}^m a_{ij}\right)^2.
\end{array}
\]
Hence
\[
\begin{array}{lll}
\frac{\dim y^\LC}{\dim y^{\GC}} \le \frac{1}{m}\left(1+\frac{1}{n-1}\right) & \Leftrightarrow & m\dim y^\LC \le \left(1+\frac{1}{n-1}\right)\dim y^\GC \\
& \Leftrightarrow & m^2k^2+2m-m\sum_{j=1}^r\sum_{i=1}^m a_{ij}^2 \le \\
&&  m^2k^2 -  \sum_{j=1}^r\left(\sum_{i=1}^m a_{ij}\right)^2 + \frac{1}{n-1}\dim y^\GC \\
& \Leftrightarrow & 2m \le \sum_j \left(m\sum_{i=1}^m a_{ij}^2 - \left(\sum_{i=1}^m a_{ij}\right)^2\right)
+ \frac{\dim y^\GC}{n-1}.
\end{array}
\]
The right hand side of the above inequality is equal to a sum of squares of the form $(a_{ij}-a_{i'j})^2$, and this sum is at least  $2m-2$. Hence the above inequality holds if $\frac{1}{n-1}\dim y^\GC \ge 2$, which is true since the minimal dimension of a nontrivial semisimple class in $\GC$ is $2n-2$. By Lemma \ref{sscomp} this proves (\ref{alin}). \hal

\vspace{6mm}
For the proof of Proposition \ref{fominclass} we need the following two lemmas.

\begin{lem}\label{albd}
Let $\KK$ be an algebraically closed field, not of characteristic $2$, let $n\ge 2$ and let $\GC$ be one of the simple algebraic groups $Sp_{2n}(\KK)$, $SO_{2n+1}(\KK)$, $SO_{2n+2}(\KK)$. Suppose $n = kr$, and let $\LC = GL_k(\KK)^r \times T$ be a Levi subgroup of $\GC$, where $T$ is a torus of rank $0,0,1$, respectively. Then $\a(\LC) \le \frac{1}{2r}$.
\end{lem}

\pf
For notational convenience we deal with the case where $\GC = Sp_{2n}(\KK)$. The other cases are entirely similar.

Write $\LC = \prod_{i=1}^r GL(V_i)$, where $\dim V_i = k$ for all $i$.
Let $u \in \LC$ be a unipotent element, and for $1\le l \le r$ write
\[
u^{V_l} = \sum_i J_i^{a_{il}},
\]
so that $\sum_i ia_{il}=k$ for $1\le l\le r$. Then as an element of $\GC = Sp_{2n}(\KK)$, we have
\[
u = \sum_i J_i^{2(a_{i1}+\cdots +a_{ir})}.
\]
To establish the lemma we need to show that $\frac{\dim u^\LC}{\dim u^{\GC}} \le \frac{1}{2r}$. Observe that
\begin{equation}\label{keyinu}
\begin{array}{ll}
\frac{\dim u^\LC}{\dim u^{\GC}} \le \frac{1}{2r} & \Leftrightarrow \dim \CB_\GC(u)-2r \dim \CB_\LC(u) \le \dim \GC - 2r \dim \LC \\
  & \Leftrightarrow \dim \CB_\GC(u)-2r \dim \CB_\LC(u) \le kr.
\end{array}
\end{equation}
By \cite[3.1]{LS},
\[
\begin{array}{ll}
\dim \CB_\LC(u) = & \sum_i i(a_{i1}^2+\cdots +a_{ir}^2) + 2\sum_{i<j}i(a_{i1}a_{j1}+\cdots +a_{ir}a_{jr}), \\
\dim \CB_\GC(u) = & 2\sum_ii(a_{i1}+\cdots +a_{ir})^2+4\sum_{i<j}i(a_{i1}+\cdots +a_{ir})(a_{j1}+\cdots +a_{jr}) \\
& + \sum_{i\;odd}(a_{i1}+\cdots +a_{ir}).
\end{array}
\]
Hence
\[
\begin{array}{ll}
 \dim \CB_\GC(u)-2r \dim \CB_\LC(u) = & 2\sum_i i\left(\left(\sum_l^r a_{il}\right)^2-r\sum_l^ra_{il}^2\right) \\
                                                        & + 4\sum_{i<j} i\left(\left(\sum_l a_{il}\right)\left(\sum_l a_{jl}\right) -
r\sum_l a_{il}a_{jl}\right) \\
                                                       & + \sum_{i\;odd}\sum_l a_{il}
\end{array}
\]
This is equal to
\[
-2\sum_i\sum_{x<y}i(a_{ix}-a_{iy})^2 - 4\sum_{i<j}\sum_{x<y}i(a_{ix}-a_{iy})(a_{jx}-a_{jy})+
\sum_{i\;odd}\sum_l a_{il}.
\]
For $x<y$, write $d_{ixy}:= a_{ix}-a_{iy}$, and define
\[
e_{xy} : = \sum_id_{ixy}^2 + 2\sum_{i<j}d_{ixy}d_{jxy},
\]
so that
\begin{equation}\label{fin}
 \dim \CB_\GC(u)-2r \dim \CB_\LC(u) = -2\sum_{x<y}e_{xy}+\sum_{i\;odd}\sum_l a_{il}.
\end{equation}
From the definition of $e_{xy}$, we see that
\[
e_{xy} = \left(d_{1xy}+\cdots +d_{rxy}\right)^2 + \left(d_{2xy}+\cdots +d_{rxy}\right)^2 +\cdots +d_{rxy}^2.
\]
In particular $e_{xy}\ge 0$, and so by (\ref{fin}),
\[
 \dim \CB_\GC(u)-2r \dim \CB_\LC(u) \le \sum_{i\;odd}\sum_l a_{il} \le  kr.
\]
This proves the inequality in  (\ref{keyinu}), completing the proof of the lemma. \hal

\begin{lem}\label{otherlem}
Let $\KK$ be algebraically closed, not of characteristic $2$, and let $X = SO_{2n}(\KK) < Y = SO_{2n+1}(\KK)$. Then for any non-identity unipotent element $u \in X$,
\[
\frac{\dim u^Y}{\dim u^X} \le 1+\frac{1}{n}.
\]
\end{lem}

\pf Let $u = \sum J_i^{a_i} \in X$ with $\sum ia_i = 2n$. Using \cite[3.1]{LS}, we see that $\dim \CB_Y(u) = \dim \CB_X(u) + \sum a_i$. Hence, writing $s = \dim [V,u] = 2n-\sum a_i$, we have
\[
\frac{\dim u^Y}{\dim u^X} = \frac{\dim Y - \dim \CB_X(u) -2n+s}{\dim X-\dim \CB_X(u)} = 1+\frac{s}{\dim u^X}.
\]
Also $\dim u^X \ge ns$ by \cite[3.4]{lish99}. The conclusion follows. \hal

\vspace{6mm}

\no {\bf Proof of Proposition \ref{fominclass}}

\vspace{2mm}
Let $m \ge 3$ be odd, $n>\frac{m^3}{3}$ and $q \equiv 1 (\mod 2m)$. Write $r = \frac{m-1}{2}$ and choose
$0\le t < m$ such that
\[
2n = km+t = 2rk+k+t.
\]
 Let $\la_1^{\pm 1},\ldots, \la_r^{\pm 1}$ be the $m^{th}$ roots of 1 in $\F_q$ that are distinct from $\pm 1$.
Write $G = G_n(q) = \GC^F$, where $\GC = Sp_{2n}(\KK)$ or $SO_{2n+\nu}(\KK)$ (where $\nu = 1$ or 2), and let $V = V_{2n+\nu}(\KK)$ be the natural module for $\GC$ (where we set $\nu=0$ in the symplectic case).  With respect to a suitable standard basis, $G$ contains the following element $x$ of order $m$:
\begin{equation}\label{eltx}
x = {\rm diag}\left(\la_1 I_{k},\la_1^{-1}I_{k}, \ldots ,\la_r I_{k}, \la_r^{-1} I_{k}, I_{k+t+\nu}\right).
\end{equation}
Then $\CB_G(x) =  L = \LC^F$, where
\[
\LC = \left\{\begin{array}{l} (GL_k)^r \times Sp_{k+t}, \hbox{ if } \GC = Sp_{2n} \\
(GL_k)^r \times SO_{k+t+\nu}, \hbox{ if } \GC = SO_{2n+\nu} \end{array}\right.
\]
To establish part (ii) of the proposition, it is sufficient to prove
\begin{equation}\label{suff1}
\dim \CB_\GC(x)-\frac{1}{m}\dim \GC < \frac{m^2}{2}.
\end{equation}
We show this in the case where $ \GC = SO_{2n+2}$ and leave the other similar cases to the reader. In this case,
\[
\begin{array}{ll}
\dim \CB_\GC(x)-\frac{1}{m}\dim \GC & = rk^2+\frac{1}{2}(k+t+1)(k+t+2) - \frac{(n+1)(2n+1)}{2r+1} \\
                                                          & = \frac{r(t+1)(t+2)}{2r+1} \\
                                                          & \le \frac{rm(m+1)}{m} < \frac{m^2}{2},
\end{array}
\]
proving (\ref{suff1}).

Finally we prove part (iii), which amounts to showing
\begin{equation}\label{suff2}
\a(\LC)\le \frac{1}{2r}\left(1+\frac{2}{n}\right).
\end{equation}
Write $\LC = \LC_1\times \LC_2$, where $\LC_1 = (GL_k)^r$ and  $\LC_2 = Sp_{k+t}$ or $SO_{k+t+\nu}$.

Let $u = u_1u_2 \in \LC$ with $u_i \in \LC_i$ unipotent elements. Obviously $\dim u^\LC = \dim u_1^{\LC_1} + \dim u_2^{\LC_2}$. Now $\GC$ has a subgroup $D = Sp_{2kr}$ or $SO_{2kr}$ containing $\LC_1$ and commuting with $\LC_2$. We claim that
\begin{equation}\label{al1bd}
\frac{\dim u_1^{\LC_1}}{\dim u_1^D} \le \frac{1}{2r}\left(1+\frac{2}{n}\right).
\end{equation}
Indeed, this follows directly from Lemma \ref{albd} in the case where $D = Sp_{2kr}$; and in the orthogonal case Lemma \ref{albd} gives $\frac{\dim u_1^{\LC_1}}{\dim u_1^{SO_{2kr+1}}} \le \frac{1}{2r}$, and so using Lemma \ref{otherlem} we have
\[
\frac{\dim u_1^{\LC_1}}{\dim u_1^D} \le \frac{1}{2r}\frac{\dim u_1^{SO_{2kr+1}}}{\dim u_1^D} \le \frac{1}{2r}\left(1+\frac{1}{kr}\right).
\]
It is easy to see that $kr>\frac{n}{2}$, so (\ref{al1bd}) follows.

Now \cite[Lemmas 4.6, 4.9]{blst1} give $\dim u^\GC \ge \dim u_1^D+\dim u_2^\GC$, whence
\[
\frac{\dim u^\LC}{\dim u^\GC} \le \frac{\dim u_1^{\LC_1}+\dim u_2^{\LC_2}}{ \dim u_1^D+\dim u_2^\GC}.
\]
Hence (\ref{suff2}) follows from (\ref{al1bd}), together with the following inequality:
\begin{equation}\label{al2bd}
\frac{\dim u_2^{\LC_2}}{\dim u_2^\GC} \le \frac{1}{2r}.
\end{equation}
To prove (\ref{al2bd}), let $s = \dim [V, u_2]$, the dimension of the commutator space of $u_2$. By \cite[3.4]{lish99}, we have
\[
\dim u_2^{\LC_2} \le \frac{s}{2}(2k+2t+2\nu-s+1),\;\;
\dim u_2^\GC \ge s(2n+\nu-s),
\]
and so
\[
\frac{\dim u_2^{\LC_2}}{\dim u_2^\GC} \le \frac{2k+2t+2\nu-s+1}{4kr+2k+2t+2\nu-2s}.
\]
Write $f(s)$ for the right hand side of the above inequality. Then $f'(s)<0$ if $k>t+\nu+1$, which holds if
$2n=km+t > (m+\nu)m+m-1$ (recall that $t\le m-1$), and this does hold because of the hypothesis $n>\frac{m^3}{3}$.
Hence the maximum value of $f(s)$ is
\[
f(1) = \frac{2k+2t+2\nu}{4kr+2k+2t+2\nu-2} \le \frac{k+t+2}{2kr+k+t+1}.
\]
Now check that the assumption $n>\frac{m^3}{3}$ implies that the right hand side is less than $\frac{1}{2r}$, proving (\ref{al2bd}). This completes the proof of the proposition. \hal

\subsection{Proof of Theorem \ref{fuchsian}}

Let $\G$ be a co-compact Fuchsian group of genus $g$ having $d$ elliptic generators $g_1,\ldots ,g_d$ of orders $m_1,\ldots ,m_d$, and define the measure $\mu = \mu(\G)$ as in (\ref{meas}). Recall that we set $v=2$ if $\G$ is oriented and $v=1$ if not; also we assume $\G$ is not virtually abelian.

We will repeatedly use the next result, taken from \cite[3.2]{fuchs2}. For a finite group $G$, and
${\bf C} = (C_1, \ldots , C_d)$  a $d$-tuple of conjugacy classes $C_i$ of $G$ with representatives $c_i$, set
\[
{\rm Hom}_{\bf C}(\G,G) = \{ \phi \in {\rm Hom}(\G,G): \phi(g_i) \in C_i
\hbox{ for } i = 1, \ldots ,d \}.
\]

\begin{lem}\label{count} {\rm (i)} If $\G$ is oriented, then
\[
|{\rm Hom}_{\bf C}(\G,G)| = |G|^{2g-1} |C_1|\cdots |C_d|
\sum_{\chi \in \Irr(G)} {{\chi(g_1)\cdots \chi(g_d)} \over {\chi(1)^{d-2+2g}}}.
\]
{\rm (ii)} If $\G$ is non-oriented, then
\[
|{\rm Hom}_{\bf C}(\G,G)| = |G|^{g-1} |C_1|\cdots |C_d|
\sum_{\chi \in \Irr(G)} \i(\c)^g {{\chi(g_1)\cdots \chi(g_d)} \over {\chi(1)^{d-2+g}}},
\]
where $\i(\c) \in \{0,1,-1\}$ is the Schur indicator of $\c$.
\end{lem}

We will also need the following extension of Theorem \ref{zetabd}.

\begin{lem}\label{zetagl} For real $s>0$, define
\[
\zeta_0^{GL_n(q)}(s) = \sum_{\c \in \Irr(GL_n(q)),\c(1)>1} \c(1)^{-s}.
\]
Then for $s > {2\over {n-1}}$,
we have $\zeta_0^{GL_n(q)}(s) \go 0$ as $q \go \infty$.
\end{lem}

\pf First note that given an irreducible character $\c$ of $SL_n(q)$, the
induced character $\c\uparrow GL_n(q)$ has at most $q-1$ irreducible constituents,
each of degree at least $\c(1)$; moreover, every irreducible character of
$GL_n(q)$ occurs as such a constituent. Recall that
\[
\zeta^{SL_n(q)}(s) = \sum_{\c\in \Irr(SL_n(q))} \c(1)^{-s}.
\]
It follows that $\zeta_0^{GL_n(q)}(s) \le (q-1)\cdot
(\zeta^{SL_n(q)}(s)-1)$.
Therefore it suffices to show that for $s > {2\over {n-1}}$, we have
$\zeta^{SL_n(q)}(s) = 1+o(q^{-1})$.
The function $\zeta^{SL_n(q)}(s)$ is investigated in
\cite[Section 2]{chardeg}, and we make
the following slight adjustments to the argument
given there to complete the proof. Adopting the notation of
\cite[Section 2]{chardeg}, we have
\[
\zeta^{SL_n(q)}(s) = \sum_{k=0}^{n-1} \D_k(t),
\]
where
\[
\D_k(t) = \sum_{\c \in \EC_k} \c(1)^{-s}
\]
(the set $\EC_k$ is defined in \cite[p.68]{chardeg}).

Each nontrivial irreducible character of $SL_n(q)$ has degree
at least $cq^{n-1}$, and hence $\D_0(s) = 1+O(q^{-(n-1)s})$, which is $1+o(q^{-1})$
for $s> {2\over {n-1}}$.

For $k \ge 1$,
\cite[Lemma 2.4]{chardeg} shows that
\[
\D_k(s) \le \hbox{max}_{\sum k_i = n-1-k} cq^{k-{1\over 2}s(n(n-1)-
\sum k_i(k_i+1))}.
\]
Hence it suffices to show that for $s > {2\over {n-1}}$
and $\sum k_i = n-1-k$, we have $k-{1\over 2}s(n(n-1)- \sum k_i(k_i+1))<-1$.
Clearly $\sum k_i(k_i+1) \le (n-k)(n-k-1)$, so it is enough to have
$2(k+1) < s(n(n-1) - (n-k)(n-k-1)$, which is true when $s>{2\over {n-1}}$.
\hal

\vspace{6mm}
\no {\bf Proof of Theorem \ref{fuchsian}, part (i) }

 Define $N_1(\G)$ as in Theorem \ref{fuchsian}. Suppose $n \ge N_1(\G)$ and $q \equiv 1 (\mod m_i)$ for all $i$, and let $G = SL_n(q)$. For each $i$, let $x_i\in G$ be an element of order $m_i$ defined as in Proposition \ref{fominsl}. Then
\begin{equation}\label{chineq}
|\c(x_i)| < f(n)\,\chi(1)^{\frac{1}{m_i}\left(1+\frac{1}{n-1}\right)} \hbox{ for all }\c \in \Irr(G).
\end{equation}
We claim also that if $C_i = x_i^G$, then
\begin{equation}\label{classineq}
|C_i| > q^{n^2(1-\frac{1}{m_i})-m_i}.
\end{equation}
To prove this, consider the centralizer of $x_i$, as given in Proposition \ref{fominsl}(ii). If $n=km_i$ with $m_i$ even and $k$ odd, then
\[
|C_i| = \frac{GL_{km_i}(q)}{GL_{k+1}(q) \times GL_{k-1}(q)\times GL_k(q)^{m_i-2}} > q^{k^2m_i^2-k^2m_i-2},
\]
which gives (\ref{classineq}). Otherwise, $n=km_i+s$ with $0\le s<m_i$, and
\[
|C_i| = \frac{GL_{km_i+s}(q)}{GL_{k+1}(q)^s \times GL_k(q)^{m_i-s}} > q^{n^2(1-\frac{1}{m_i})-s(1-\frac{s}{m_i})},
\]
which implies (\ref{classineq}).

In order to make use of Lemma \ref{count}, define
\[
\Sigma = \sum_{\c \in \Irr(G), \c \ne 1} \frac{|\chi(x_1)\cdots \chi(x_d)|}{\chi(1)^{d-2+vg}}.
\]
By (\ref{chineq}), we have $\Sigma \le f(n)^d\sum_{\c\ne 1} \c(1)^{-\g} = f(n)^d (\zeta^G(\g)-1)$, where $\zeta^G$ is as defined in (\ref{zetadef}) and
\[
\g = d-2+vg-\sum \frac{1}{m_i}-\frac{1}{n-1}\sum \frac{1}{m_i} = \mu -\frac{1}{n-1}\sum \frac{1}{m_i}.
\]
By Theorem \ref{zetabd}, $\zeta^G(\g) \rightarrow 1$ as $q \rightarrow \infty$ provided $\g > \frac{2}{n}$, and this inequality holds by our assumption that $n \ge N_1(\G)$. Hence for sufficiently large $q$ we have $\Sigma \le \frac{1}{2}$, and so by Lemma \ref{count},
\[
|{\rm Hom}_{\bf C}(\G,G)| \ge \frac{1}{2}|G|^{vg-1} |C_1|\cdots |C_d|.
\]
Using (\ref{classineq}), it follows that there is a positive constant $c$ such that
$|{\rm Hom}_{\bf C}(\G,G)| \ge cq^{t}$, where
\[
t = (n^2-1)\left(vg-1+\sum \left(1-\frac{1}{m_i}\right)\right) + \sum \left(1-\frac{1}{m_i}\right) - \sum m_i.
\]
Then $t = (n^2-1)(\mu+1) + \sum \left(1-\frac{1}{m_i}\right) - \sum m_i$, and the conclusion of Theorem \ref{fuchsian}(i) follows.

\vspace{6mm}
\no {\bf Proof of Theorem \ref{fuchsian}, part (ii) }

 As usual let $\GC = GL_n(\KK)$ and $G = GL_n(q) = \GC^F$. Let $m\ge 2$ and suppose $q \equiv 1 (\mod m)$ and $n\ge m$. The centralizers in $G$ of elements of order $m$ are Levi subgroups $L = \LC^F$. Let $L_0 = \LC_0^F$ be such a Levi subgroup of minimal order.

Let $y \in G$ be an element of order $m$, and $L = \LC^F = \CB_G(y)$. Then Theorems \ref{charbd} and \ref{ratio} show that for any $\c \in \Irr(G)$,
\begin{equation}\label{prod1}
|y^G|\,|\c(y)| \le \frac{|G|}{|L|}f(n)\c(1)^{\frac{1}{2}\left(1+\frac{\dim \LC}{\dim \GC}\right)}.
\end{equation}

We now assert that if $\dim \LC > \dim \LC_0$, then for large $q$,
\begin{equation}\label{prod2}
\frac{\c(1)^{\frac{1}{2}\left(1+\frac{\dim \LC}{\dim \GC}\right)}}{|L|} \le
\frac{\c(1)^{\frac{1}{2}\left(1+\frac{\dim \LC_0}{\dim \GC}\right)}}{|L_0|}.
\end{equation}
To see this, note that obviously $\c(1) < q^{\frac{1}{2}\dim \GC}$. This implies that
\[
\frac{\c(1)^{\frac{\dim \LC_0}{\dim \GC}}}{|L_0|} \frac{|L|} {\c(1)^{\frac{\dim \LC}{\dim \GC}}}
\sim \frac{q^{\dim \LC-\dim \LC_0}}{\c(1)^{\frac{\dim \LC-\dim\LC_0}{\dim \GC}}} \ge
 \frac{q^{\dim \LC-\dim \LC_0}}{q^{\frac{1}{2}(\dim \LC-\dim \LC_0)}} \ge  q^{\frac{1}{2}},
\]
which establishes (\ref{prod2}).

Next we claim that
\begin{equation}\label{prod3}
\frac{n^2}{m} \le \dim \LC_0 \le \frac{n^2}{m}+\frac{m}{4}.
\end{equation}
To derive the lower bound, observe that since $\LC_0$ is the centralizer of an element of order $m$, we have $\LC_0 = \prod_{i=1}^r GL_{n_i}(\KK)$, where $\sum n_i = n$ and $r \le m$. Hence
$$m \dim \LC_0 = m\sum_{i=1}^r n_i^2 \ge r\sum_{i=1}^r n_i^2 \ge \left(\sum n_i\right)^2 = n^2,$$
proving the lower bound. For the upper bound, write
$n=km+s$ with $0\le s<m$, and observe that $L = \LC^F = GL_{k+1}(q)^s \times GL_k(q)^{m-s}$ is the centralizer of an element of order $m$, as in Proposition \ref{fomin}. Hence
\[
\begin{array}{ll}
\dim \LC_0 \le \dim \LC & = mk^2+2ks+s \\
                                     & = \frac{n^2}{m}+\frac{ms-s^2}{m} \\
                                     & \le \frac{n^2}{m}+\frac{m}{4},
\end{array}
\]
giving the upper bound in (\ref{prod3}).

For $1\le i\le d$, let $y_i\in G$ be an element of order $m_i$, and let $C_i = y_i^G$. Applying (\ref{prod1})--(\ref{prod3}), we see that for any $\c \in \Irr(G)$ with $\c(1)>1$,
\[
|C_i|\,|\c(y_i)| \le \frac{f(n)|G| \c(1)^{\frac{1}{2}\left(1+\frac{1}{m_i}+\frac{m_i}{4n^2}\right)}}{q^{\frac{n^2}{m_i}}}
\le  f(n) q^{n^2\left(1-\frac{1}{m_i}\right)} \c(1)^{\frac{1}{2}\left(1+\frac{1}{m_i}+\frac{m_i}{4n^2}\right)}.
\]
Hence
\[
\prod_{i=1}^d|C_i|\sum_{\c(1)>1} \frac{|\c(y_1)\cdots \c(y_d)|}{\c(1)^{d-2+vg}} \le
f(n)^dq^{n^2\sum \left(1-\frac{1}{m_i}\right)} \sum_{\c(1)>1}\c(1)^{-\g},
\]
where $\g = d-2+vg - \frac{1}{2}\sum \left(1+\frac{1}{m_i}+\frac{m_i}{4n^2}\right)$. By Lemma \ref{zetagl},
$$\sum_{\c(1)>1}\c(1)^{-\g} = \zeta_0^G(\g) \go 0$$
as $q \go \infty$ provided $\g > \frac{2}{n-1}$. Now
\[
2\g = \mu-2+vg - \sum \frac{m_i}{4n^2} \ge \mu-2- \frac{d}{4n},
\]
from which we see that  if $\mu > 2$ and $n-1 \ge \frac{d+16}{4(\mu-2)}$, then $\g > \frac{2}{n-1}$.
So under these conditions (which are in the hypothesis of Theorem \ref{fuchsian}(ii)), for $q$ sufficiently large we have
\[
\prod_{i=1}^d|C_i|\sum_{\c(1)>1} \frac{|\c(y_1)\cdots \c(y_d)|}{\c(1)^{d-2+vg}} \le
\frac{1}{2}f(n)^dq^{n^2\sum \left(1-\frac{1}{m_i}\right)}.
\]
Adding in the contribution from the $q-1$ linear characters of $G$, it now follows from Lemma \ref{count} that
\[
|{\rm Hom}_{\bf C}(\G,G)| \le  |G|^{vg-1}\,\left((q-1)\prod |C_i| +
\frac{1}{2}f(n)^d q^{n^2\sum \left(1-\frac{1}{m_i}\right)}\right).
\]
From (\ref{prod3}) we have $|C_i| \le q^{n^2\left(1-\frac{1}{m_i}\right)}$, and hence for large $q$,
\begin{equation}\label{prod4}
|{\rm Hom}_{\bf C}(\G,G)| \le (q-1)\,q^{n^2\left(vg-1+\sum \left(1-\frac{1}{m_i}\right)\right)} = (q-1)\,q^{n^2(\mu+1)}.
\end{equation}
To complete the proof of Theorem \ref{fuchsian}(ii), observe that $|{\rm Hom}(\G,G)|$ is equal to the sum
$\sum_{\bf C} |{\rm Hom}_{\bf C}(\G,G)|$ over all $d$-tuples ${\bf C}$ of classes of elements of orders $m_1,\ldots ,m_d$. The number of such classes is bounded by a function of $n$, so the conclusion of Theorem \ref{fuchsian}(ii) follows
from (\ref{prod4}). \hal

\subsection{Proof of Theorem \ref{fuchsclass}}

First we prove part (i) of the theorem. This runs along similar lines to the proof of Theorem \ref{fuchsian}(i).
As in the hypothesis of Theorem \ref{fuchsclass}, let $\G$ be a Fuchsian group with $m_i$ odd for all $i$, and $\mu > t(\G)$. Suppose $n > N_3(\G)$ and $q\equiv 1  (\mod 2m_i)$ for all $i$, and let $G = G_n(q)$ be one of the groups
$Sp_{2n}(q)$, $\O_{2n+1}(q)$, $\O_{2n+2}^\pm (q)$.

For each $i$, let $x_i\in G$ be an element of order $m_i$, defined as in Proposition \ref{fominclass}, and let $C_i = x_i^G$. Then for all $\c \in \Irr(G)$, Proposition \ref{fominclass}(iii) gives
\[
|\c(x_i)| < f(n)\,\chi(1)^{\frac{1}{m_i-1}\left(1+\frac{2}{n}\right)}.
\]
Hence
\[
 \sum_{\c \in \Irr(G), \c \ne 1} \frac{|\chi(x_1)\cdots \chi(x_d)|}{\chi(1)^{d-2+vg}} \le f(n)^d\sum_{\c\ne 1} \c(1)^{-\g},
\]
where $\g = d-2+vg-\left(1+\frac{2}{n}\right)\sum \frac{1}{m_i-1} = \mu-t(\G)-\frac{2}{n}\sum \frac{1}{m_i-1}$. Since $n>N_3(\G)$ by hypothesis, we have $\g > \frac{1}{n} = \frac{2}{h}$ (where $h$ is the Coxeter number of $G$), so by Theorem  \ref{zetabd}, for large $q$ the above sum is less than
$\frac{1}{2}$.

Proposition \ref{fominclass}(ii) gives $|C_i| > c|G|^{1-\frac{1}{m_i}}q^{-\frac{m_i^2}{2}}$.
It now follows using Lemma \ref{count} that
\[
\begin{array}{ll}
|{\rm Hom}_{\bf C}(\G,G)| & \ge \frac{1}{2}|G|^{vg-1} |C_1|\cdots |C_d| \\
                                            & \ge c'|G|^{vg-1+\sum\left(1-\frac{1}{m_i}\right)}q^{-\frac{1}{2}\sum m_i^2} \\
                                             & = c'|G|^{1+\mu}q^{-\frac{1}{2}\sum m_i^2},
\end{array}
\]
where $c'$ is a positive absolute constant. This completes the proof.

\vspace{4mm}
Now we prove part (ii) of Theorem \ref{fuchsclass}.
As in the hypothesis, assume that all $m_i$ are odd and that $\mu > 2$ and $n> N_4(\G)$. Let $G = G_n(q) = \GC^F$, where $\GC = G_n(\KK)$ and $\KK$ is algebraically closed. Suppose $m\ge 3$ is an odd number such that $q \equiv 1 (\mod 2m)$ and $n \ge m$. Let $y\in G$ have order $m$. Then $\CB_G(y) = L = \LC^F$, a split Levi subgroup; let $L_0 = \LC_0^F$ be such a Levi of minimal order. As in the proof of Theorem \ref{fuchsian}(ii), the inequalities (\ref{prod1}) and (\ref{prod2}) hold. We also claim as in (\ref{prod3}) that
\begin{equation}\label{loin}
\frac{\dim \GC}{m}-1 \le \dim \LC_0 \le  \frac{\dim \GC}{m}+\frac{m^2}{2}.
\end{equation}
To see this, note first that the upper bound is implied by the centralizer in Proposition \ref{fominclass}. We sketch a proof of the lower bound for the case $\GC = Sp_{2n}$:  for a Levi $\LC = \CB_\GC(y)$ as above, writing $r=\frac{m-1}{2}$ we have
\[
\LC = Sp_{2t} \times \prod_{i=1}^r GL_{k_i},
\]
where $n = t+\sum_{i=1}^r k_i$. Hence
\[
\begin{array}{ll}
\dim \LC & = 2t^2+t + \sum k_i^2 \\
              & \ge 2t^2+t+\frac{1}{r}(n-t)^2 \\
              & = \frac{1}{r}\left((2r+1)t^2-(2n-r)t+n^2\right),
\end{array}
\]
and it is straightforward to see that the right hand side is at least $\frac{2n^2+n}{m}-1$, proving the lower bound in (\ref{loin}).

Now let $y_i\in G$ be elements of order $m_i$ for $1 \le i \le d$, and let $C_i = y_i^G$. Then by (\ref{prod1}), (\ref{prod2}) and (\ref{loin}),
\[
|C_i|\,|\c(y_i)| \le \frac{f(n)|G| \c(1)^{\frac{1}{2}\left(1+\frac{1}{m_i}+\frac{m_i^2}{2\dim \GC}\right)}}
{q^{-1}|G|^{\frac{1}{m_i}}}.
\]
Hence
\[
\prod_{i=1}^d|C_i|\sum_{\c(1)>1} \frac{|\c(y_1)\cdots \c(y_d)|}{\c(1)^{d-2+vg}} \le
f(n)^dq^d |G|^{\sum \left(1-\frac{1}{m_i}\right)} \sum_{\c \ne 1}\c(1)^{-\g},
\]
where $\g = d-2+vg-\frac{1}{2}\sum\left(1+\frac{1}{m_i}+\frac{m_i^2}{4n^2}\right)$. Then
\[
2\g = \mu +vg-2-\sum \frac{m_i^2}{4n^2} \ge \mu-2-\frac{d}{4n}
\]
(recall that $n>N_4(\G) \ge m_i^2$).
Hence the assumption $n>N_4(\G)\ge \frac{d+8}{4(\mu-2)}$ implies that
$\g > \frac{1}{n} = \frac{2}{h}$. Now Theorem \ref{zetabd} shows that for large $q$, and a suitable function $g(n)$,
\[
|{\rm Hom}_{\bf C}(\G,G)| \le g(n)q^d |G|^{vg-1+\sum \left(1-\frac{1}{m_i}\right)} = g(n)q^d|G|^{\mu+1}.
\]
This completes the proof of Theorem \ref{fuchsclass}. \hal

\subsection{Deduction of Corollaries \ref{dimhomgl} and \ref{dimhomclass}}

It is routine to deduce these corollaries from Theorems \ref{fuchsian} and \ref{fuchsclass}. We do this for Corollary \ref{dimhomgl} and leave the other very similar deduction to the reader.

Let $\G, n$ and $\KK$ be as in Corollary \ref{dimhomgl}. First suppose $\KK = \bar \F_p$, where $p$ is a prime. Let $V$ be the variety $R_{n,\KK}(\G)$. For a power $q$ of $p$, let $V(q)$ denote the set of $\F_q$-rational points in $V$. Let $f = \dim V$, let $V_i\,(1\le i\le r)$ be the set of irreducible components of $V$ of dimension $f$, and $W_i\,(1\le i\le s)$ the irreducible components of dimension less than $f$.

Choose a power $q_0$ of $p$ such that all $V_i$ and $W_i$ are defined over $\F_{q_0}$. Then for $q = q_0^k$ the Lang-Weil estimate \cite{LW} gives $|V_i(q)| = q^f+O(q^{f-\frac{1}{2}})$ for $1\le i\le r$, while
$|W_i(q)| = O(q^{f-1})$ for $1\le i\le s$. Since there are infintely many powers $q_0^k$ that are congruent to 1 modulo $m_i$ for all $i$, the parts (i) and (ii) of Corollary \ref{dimhomgl} now follow from Theorem \ref{fuchsian}, in the case where $\KK = \bar \F_p$. It follows that these parts also hold for all algebraically closed fields of characteristic $p$. They follow in characteristic zero as well, since it is well known that the dimension of a variety in characteristic zero coincides with the dimension of its reduction modulo $p$ for all large primes $p$.

Finally, part (iii) of Corollary \ref{dimhomgl} follows from part (ii), together with the fact that, in the above notation, $|V(q)| = O(q^f)$, which again follows using \cite{LW}.

\subsection{Proof of Theorem \ref{excep1}}\label{lines}

Again let $\G$ be a co-compact, non-virtually abelian Fuchsian group of genus $g$ having $d$ elliptic generators
$g_1,\ldots ,g_d$ of orders $m_1,\ldots ,m_d$, and define the measure $\mu = \mu(\G)$ as in (\ref{meas}). Assume that $m_1\cdots m_d$ is coprime to 30.

Let $\GC$ be an adjoint simple algebraic group of exceptional type over an algebraically closed field $\KK$ of good characteristic $p$ not dividing $m_1\cdots m_d$. Let $q$ be a power of $p$ such that $q \equiv 1 (\mod m_i)$ for all $i$, and let $F$ be a Frobenius endomorphism of $\GC$ such that $G = G(q) = \GC^F$ is a finite group of Lie type over $\F_q$. For each $i$ define $J_{m_i}(\GC) = \{x \in \GC : x^{m_i}=1\}$, and let $j_{m_i}:=j_{m_i}(\GC)
= \dim J_{m_i}(\GC)$.

We will prove that there are positive constants $c_1,c_2$ such that for sufficiently large $q$,
\begin{equation}\label{ulbds}
c_1|G|^{vg-1}q^{\sum j_{m_i}} \le |{\rm Hom}(\G,G)| \le c_2|G|^{vg-1}q^{\sum j_{m_i}}.
\end{equation}
Just as in the previous subsection, this implies Theorem \ref{excep1}.

For the proof of (\ref{ulbds}), the following bounds on the dimensions $j_{m_i}$ will be useful; they are taken from \cite[Theorem 1]{law}, bearing in mind that each $m_i \ge 7$ by assumption:
\begin{equation}\label{jbds}
\begin{array}{|c|ccccc|}
\hline
\GC & E_8 & E_7& E_6& F_4& G_2 \\
\hline
j_{m_i}\ge & 212 & 114 & 66 & 44 & 12 \\
\hline
\frac{2}{h} & \frac{1}{15} & \frac{1}{9} & \frac{1}{6}& \frac{1}{6}& \frac{1}{3} \\
\hline
\end{array}
\end{equation}
For the reader's convenience we have also included the values of $\frac{2}{h}$ where $h$ is the Coxeter number of $\GC$.

First we establish the lower bound on $|{\rm Hom}(\G,G)|$ in (\ref{ulbds}). For each $i$, since $q \equiv 1 (\mod m_i)$, there exists $x_i \in G$ such that $\dim x_i^\GC = j_{m_i}$ and $\CB_G(x_i) = L_i = \LC_i^F$, where $\LC_i$ is an $F$-stable split Levi subgroup in $\GC$. The possible Levi subgroups $\LC_i$ (realising the maximal dimension $j_{m_i}$) are given in the tables on pp.240-241 of \cite{law}, and it follows using Theorem \ref{alphaexcep} that $\a(\LC_i) \le \frac{1}{6}$, and also $\a(\LC_i) = 0$ for $\GC = G_2$.
Hence by Theorem \ref{charbd}, provided $q$ is sufficently large we have, for $\GC \ne G_2$,
\[
\sum_{\c \in \Irr(G),\c \ne 1} \frac{|\c(x_1)\cdots \c(x_d)|}{\c(1)^{d-2+vg}} \le \sum_{\c \in \Irr(G),\c \ne 1}
\c(1)^{-(d-2+vg)+\frac{d}{6}},
\]
and the same bound without the $\frac{d}{6}$ term for $\GC=G_2$. It is easy to check that $\frac{5}{6}d-2+vg > \frac{2}{h}$ for $\GC \ne G_2$ and $d-2+vg > \frac{2}{h}$ for $G_2$ (recall that $vg+d\ge 3$ as $\G$ is not virtually abelian). Hence Theorem  \ref{zetabd} implies that the above sum tends to 0 as $q\rightarrow \infty$. It now follows from Lemma \ref{count} that for large $q$,
\begin{equation}\label{dec}
|{\rm Hom}(\G,G)| \ge \frac{1}{2}|G|^{vg-1}\prod_{i=1}^d |x_i^G|,
\end{equation}
and the lower bound in (\ref{ulbds}) follows.

Now we prove the upper bound on $|{\rm Hom}(\G,G)|$ in (\ref{ulbds}).
For each $i$ let $y_i \in G$ be an element of order $m_i$. Since $q \equiv 1 (\mod m_i)$ and $m_i$ is coprime to 30, hence is not divisible by a bad prime for $\GC$, it follows that $\CB_G(y_i) = L_i = \LC_i^F$, where $\LC_i$ is an $F$-stable split Levi subgroup of $\GC$ (see \cite[4.2.2]{GLS}).
 Let $C_i = y_i^G$, ${\mathcal C}_i = y_i^\GC$
 and assume that $\dim {\mathcal C}_1 \le \cdots \le \dim {\mathcal C}_d$. Write ${\bf C} = (C_1,\ldots, C_d)$.

From the presentation of $\G$ it is clear that $|{\rm Hom}_{\bf C}(\G,G)| \le |G|^{vg}\prod_{i=1}^{d-1}|C_i|$, and hence, provided $d\ge 2$, there is a positive constant $c$ such that
\begin{equation}\label{hombd}
|{\rm Hom}_{\bf C}(\G,G)| \le c|G|^{vg-1}q^{\sum j_{m_i}}\cdot q^{\dim \GC + \dim {\mathcal C}_1-j_{m_1}-j_{m_d}}.
\end{equation}
Also by Lemma \ref{count} together with Theorem \ref{charbd},
\[
\begin{array}{ll}
|{\rm Hom}_{\bf C}(\G,G)| &  \le |G|^{vg-1} |C_1|\cdots |C_d|
\sum_{\chi \in \Irr(G)}  \frac{|\chi(y_1)\cdots \chi(y_d)|}{\chi(1)^{d-2+vg}} \\
     & \le c'|G|^{vg-1}q^{\sum j_{m_i}} \sum_{\chi \in \Irr(G)} \c(1)^{\sum \a(\LC_i)-(d-2+vg)}.
\end{array}
\]
If $(d-2+vg)-\sum \a(\LC_i) > \frac{2}{h}$ where $h$ is the Coxeter number of $\GC$, then the upper bound in (\ref{ulbds}) follows by Theorem \ref{zetabd}. Hence we may assume that
\begin{equation}\label{zineq}
d-2+vg-\sum_{i=1}^d \a(\LC_i) \le \frac{2}{h}.
\end{equation}

Recall our assumption that $\G$ is not virtually abelian, hence $vg+d \ge 3$. If $d=0$ then $vg\ge 3$ and
(\ref{zineq}) fails. And if $d = 1$ then $vg\ge 2$ and (\ref{zineq}) gives $\a(\LC_1)\ge 1-\frac{2}{h}$, which is not possible by Theorem \ref{alphaexcep}.

Now suppose $d=2$. Then $vg\ge 1$ and (\ref{zineq}) gives
\[
\a(\LC_1)+\a(\LC_2) \ge 1-\frac{2}{h}.
\]
So we can assume that $\a(\LC_1) \ge \frac{1}{2}-\frac{1}{h}$. By Theorem \ref{alphaexcep}, the possibilities for $\LC_1$ and consequent upper bounds for $\dim {\mathcal C}_1 = \dim \GC -\dim \LC_1$ are as follows:
\[
\begin{array}{|c|c|c|c|c|c|}
\hline
\GC & E_8& E_7& E_6& F_4& G_2 \\ \hline
\LC_1' & E_7 & E_6,D_6 & D_5,A_5,D_4 & B_3,C_3 & A_1 \\ \hline
\dim {\mathcal C}_1 \le & 114 & 66 & 48 & 30 & 10\\
\hline
\end{array}
\]
In all cases we check using (\ref{jbds}) that $\dim \GC + \dim {\mathcal C}_1-j_{m_1}-j_{m_d} \le 0$. By (\ref{hombd}), this gives the upper bound in (\ref{ulbds}).

Finally, suppose that $d\ge 3$. Then (\ref{zineq}) gives
\[
\sum_{i=1}^d \a(\LC_i) \ge d-2-\frac{2}{h}.
\]
So we can assume that $\a(\LC_1) \ge 1-\frac{1}{d}\left(2+\frac{2}{h}\right) \ge \frac{1}{3}- \frac{2}{3h}$.
By Theorem \ref{alphaexcep}, the possibilities for $\LC_1$ are as follows:
\[
\begin{array}{|c|c|c|c|c|c|}
\hline
\GC & E_8& E_7& E_6& F_4& G_2 \\ \hline
&&&&&\\
\vspace{-8.5mm}\\
\LC_1\ \triangleright & E_7,D_7,E_6 & E_6,D_6,D_5,A_5,A_5' & D_5,A_5,D_4,A_4 & B_3,C_3 & A_1,\tilde A_1 \\ \hline
\dim {\mathcal C}_1 \le & 168 & 94 & 50 & 30 & 10\\
\hline
\end{array}
\]
Again, we check using (\ref{jbds}) that $\dim \GC + \dim {\mathcal C}_1-j_{m_1}-j_{m_d} \le 0$. By (\ref{hombd}), this gives the upper bound in (\ref{ulbds}).

The upper bound in (\ref{ulbds}) is now proved, and the proof of Theorem \ref{excep1} is complete. \hal

\section{Probabilistic generation of classical groups}
\label{slfuchs}

In  this section we prove Theorems \ref{fuchsgensl} and  \ref{fuchsgenclass}.

Let $\G$ be a co-compact Fuchsian group of genus $g$ having $d$ elliptic generators $g_1,\ldots ,g_d$ of orders $m_1,\ldots ,m_d$, and define the measure $\mu = \mu(\G)$ as in (\ref{meas}). Recall that we set $v=2$ if $\G$ is oriented and $v=1$ if not; also we assume $\G$ is not virtually abelian, so that $vg+d\ge 3$.

If $vg\ge 3$ then Theorems \ref{fuchsgensl}  and \ref{fuchsgenclass} follow from \cite[1.6]{fuchs2}, so we assume throughout that
\[
vg \le 2.
\]

\subsection{Proof of Theorem \ref{fuchsgensl}}

We assume the following bounds on $\mu$ and $n$, as in the hypotheses of Theorem \ref{fuchsgensl}:
\begin{equation}\label{hyps}
\mu > \nu:= {\rm max}\left(2,\,1+\sum\frac{1}{m_i}\right),\;\;
n \ge \nu N_2(\G)+2\sum m_i.
\end{equation}

Let $G = SL_n(q)$ with $q \equiv 1  (\mod m_i)$ for all $i$. By Theorem \ref{fuchsian}(i), for large $q$ we have
\[
|{\rm Hom}(\G,G)| > |G|^{\mu+1}q^{-\sum m_i}.
\]
Write $D$ for the degree of the rational function in $q$ given by the right hand side, so that
\[
D = (n^2-1)(\mu+1)-\sum m_i.
\]
Now $P_\G(G)$ is the probability that a randomly chosen element of ${\rm Hom}(\G,G)$ is an epimorphism, so clearly
\[
1-P_\G(G) \le \frac{\sum_{M < G {\mbox{\tiny{ is maximal}}}}|{\rm Hom}(\G,M)|}{|{\rm Hom}(\G,G)|}.
\]
Hence the conclusion of Theorem \ref{fuchsgensl} will follow if we show that
\begin{equation}\label{keyin}
\sum_{M < G {\mbox{\tiny{ is maximal}}}}|{\rm Hom}(\G,M)| < cq^{D-\frac{1}{2}},
\end{equation}
where $c$ is a constant. In the rest of the proof we aim to establish (\ref{keyin}).

Define ${\mathcal M}_1$ to be a set of conjugacy class representatives of the maximal subgroups $M$ of $G$ such that
\[
|M| < |G|^{\frac{\mu}{vg+d-2}}q^{-\sum m_i-1},
\]
and let ${\mathcal M}_2$ to be a set of conjugacy class representatives for the remaining maximal subgroups. For $i=1,2$ let
\[
\Sigma_i = \sum_{M \in {\mathcal M}_i} |{\rm Hom}(\G,M)| \cdot |G:M|,
\]
so that $\sum_{M\,max\,G}|{\rm Hom}(\G,M)| = \Sigma_1+\Sigma_2$.

Consider first $\Sigma_1$. For $M \in {\mathcal M}_1$,
\[
\begin{array}{ll}
|{\rm Hom}(\G,M)|\cdot |G:M| & \le |M|^{vg+d-1}|G:M| =   |M|^{vg+d-2}|G| \\
                                               & \le |G|^{\mu+1}q^{(-\sum m_i-1)(vg+d-2)} \\
                                               & \le cq^{D-1}.
\end{array}
\]
Also $|{\mathcal M}_1| \le c(n)\log \log q$ by \cite[1.3]{LMS}. It follows that for large $q$,
\begin{equation}\label{s1bd}
\Sigma_1 \le q^{D-\frac{1}{2}}.
\end{equation}

Bounding $\Sigma_2$ takes more work. Let $M \in {\mathcal M}_2$, so that
\[
|M| \ge |G|^{\frac{\mu}{vg+d-2}}q^{-\sum m_i-1}.
\]
Write $w = vg+d-2$ and $b = \sum m_i+2$. The assumption in (\ref{hyps}) that $\mu > 1+\sum\frac{1}{m_i}$ implies
that  $\mu \ge \frac{1}{2}(vg+d-1)$, which yields
$\frac{\mu}{vg+d-2} \ge \frac{1}{2}(1+\frac{1}{w})$, and hence
\[
 |G|^{\frac{\mu}{vg+d-2}}q^{-\sum m_i-1} \ge cq^{\frac{1}{2}(n^2-1)(1+\frac{1}{w})-b+1}
\ge cq^{\frac{1}{2}n^2(1+\frac{1}{w})-b}
\]
where $c$ is a positive constant.
Now $w \le d<\frac{b}{2}$ (recall we are assuming $vg\le 2$), and $n>2b$ by (\ref{hyps}), so
$\frac{1}{2}n^2(1+\frac{1}{w})-b > \frac{1}{2}(n^2+n)+1$. It follows that
\begin{equation}\label{spbd}
|M| \ge cq^{\frac{1}{2}(n^2+n)+1}.
\end{equation}
Observe also that the assumptions $vg\le 2$ and $\mu>2$ imply that $d\ge 3$, from which it is easy to see from (\ref{hyps}) that $n\ge18$.

According to the well-known theorem of Aschbacher \cite{asch}, the maximal subgroups of $G$ fall into eight classes
${\mathcal C}_i$ ($1\le i\le 8$) of ``geometric" subgroups, together with a class ${\mathcal S}$ consisting of subgroups that are almost simple modulo scalars, and act  absolutely irreducibly on $V=V_n(q)$. For large $q$ the subgroups in ${\mathcal S}$ have order less that $q^{3n}$ by \cite{Li}. Furthermore, by inspection of the subgroups in ${\mathcal C}_i$ (see \cite[Chapter 4]{KL}), the largest maximal subgroup in $\bigcup_{i=2}^8 {\mathcal C}_i$ has order at most $cq^{\frac{1}{2}(n^2+n)}$ (this is attained for the symplectic group $Sp_n(q)$ when $n$ is even). Hence it follows from (\ref{spbd}) that $M$ belongs to the class
${\mathcal C}_1$, which consists of the maximal parabolic subgroups of $G$.

So $M$ is a parabolic subgroup of $G$. Write $M = QL$, where $Q$ is the unipotent radical and $L$ a Levi subgroup. Then $L = G\cap (GL_r(q) \times GL_{n-r}(q))$ for some $r \le \frac{n}{2}$, and $|Q| = q^{r(n-r)}$. Note that
$|G:M| = |G:QL| \le c|Q|$, so
\begin{equation}\label{rad}
\begin{array}{rl}
|{\rm Hom}(\G,QL)|\cdot |G:QL|  \le & |Q|^{vg+d-1}|{\rm Hom}(\G,L)|\cdot c|Q| \\
                   = & c|Q|^{vg+d}|{\rm Hom}(\G,L)|.
\end{array}
\end{equation}

Assume first that $r>N_2(\G)$. Then by Theorem \ref{fuchsian}(ii), for $s = r$ or $n-r$ we have
\begin{equation}\label{glbds}
|{\rm Hom}(\G,GL_s(q))| < g(n) q |GL_s(q)|^{\mu+1},
\end{equation}
and so
\[
|{\rm Hom}(\G,L)| < g(n) |L|^{\mu+1}q^2q^{\mu+1}.
\]
Hence
\begin{equation}\label{yet1}
|{\rm Hom}(\G,QL)|\cdot |G:QL|  \le g(n) |Q|^{vg+d}|L|^{\mu+1}q^{\mu+3}.
\end{equation}
Write $\dim Q$, $\dim L$ for the degrees in $q$ of the polynomials $|Q|$, $|L|$, and note that
$\dim G = \dim L + 2\dim Q$. Then
\begin{equation}\label{dexp}
D = (n^2-1)(\mu+1)-\sum m_i = \left(\dim L+2\dim Q\right)(\mu+1)-\sum m_i.
\end{equation}
The degree of the right hand side of (\ref{yet1}) is
\[
(vg+d)\dim Q+(\mu+1)\dim L+\mu+3.
\]
 Now
\[
\begin{array}{ll}
& (vg+d)\dim Q+(\mu+1)\dim L+\mu+3 < D-1 \\
 \Leftrightarrow & \left(2(\mu+1)-(vg+d)\right)\dim Q > \sum m_i+\mu +4 \\
 \Leftrightarrow & \mu > \frac{1}{2}(vg+d)-1+\frac{\sum m_i+\mu+4}{2\dim Q}.
\end{array}
\]
Since we are assuming that $\mu \ge \frac{1}{2}(vg+d-1)$, the above inequality holds provided $\dim Q > \sum m_i + \mu +4$, which holds by our assumption on $n$ in (\ref{hyps}).

Hence in the case where $r>N_2(\G)$, it now follows by (\ref{yet1}) that
\begin{equation}\label{yet2}
|{\rm Hom}(\G,QL)|\cdot |G:QL| < cq^{D-1}.
\end{equation}

Now assume that $r \le N_2(\G)$. In this case (\ref{glbds}) holds only for $s=n-r$ (note that $n>2N_2(\G)$ by (\ref{hyps})), so
\[
|{\rm Hom}(\G,QL)|\cdot |G:QL| <  g(n)q |Q|^{vg+d}|GL_r(q)|^{vg+d-1}|GL_{n-r}(q)|^{\mu+1}.
\]
The right hand side has degree in $q$ at most $R$, where
\[
R:= (vg+d)\dim Q+(\mu+1)\dim L+(\sum \frac{1}{m_i})r^2+\mu +2.
\]
From the expression for $D$ in (\ref{dexp}), we have
\[
R \le D-1 \Leftrightarrow (\mu-\sum \frac{1}{m_i})\dim Q \ge (\sum \frac{1}{m_i})r^2+ \sum m_i+\mu+2.
\]
As $\dim Q = r(n-r)$, this inequality holds if
\[
n> r(1+\sum \frac{1}{m_i}) + \frac{1}{r}(\sum m_i+\mu+2),
\]
and since $r\le N_2(\G)$, this holds by (\ref{hyps}).

Hence (\ref{yet2}) holds also when $r\le N_2(\G)$. It follows that
\[
\Sigma_2 = \sum_{M \in {\mathcal M}_2} |{\rm Hom}(\G,M)| \cdot |G:M| < c'q^{D-1}.
\]
Together with (\ref{s1bd}), this proves (\ref{keyin}).

This completes the proof of Theorem \ref{fuchsgensl}. \hal

\subsection{Proof of Theorem \ref{fuchsgenclass}}

The proof runs along similar lines to the previous section, and we will omit quite a few details.
As in the hypothesis of Theorem  \ref{fuchsgenclass}, assume that
all the $m_i$ are odd, and that
$$\mu > {\rm max} \left(2,\, t(\G),\,1+\sum \frac{1}{m_i}\right).$$
Assume also that $n > N_5(\G)$ (as defined in the preamble to the theorem), and let $G = G_n(q)$ be
$Sp_{2n}(q)$, $\O_{2n+1}(q)$, or $\O_{2n+2}^\pm(q)$ with $n> N_5(\G)$ and $q \equiv 1  (\mod 2m_i)$ for all $i$.

By Theorem \ref{fuchsclass}, we have $|{\rm Hom}(\G,G)| > cq^D$, where $c$ is a positive constant and
\[
D = (\mu+1)\dim G - \frac{1}{2}\sum m_i^2
\]
(here, as in the previous section we write $\dim G$ for the degree of $|G|$ as a polynomial in $q$). Let ${\mathcal M}_1$ be a set of representatives of the conjugacy classes of maximal subgroups $M$ of $G$ such that
\[
|M| < |G|^{\frac{\mu}{vg+d-2}}q^{-\frac{1}{2}\sum m_i^2-1},
\]
and let ${\mathcal M}_2$ to be a set of conjugacy class representatives for the remaining maximal subgroups. For $i=1,2$ let
\[
\Sigma_i = \sum_{M \in {\mathcal M}_i} |{\rm Hom}(\G,M)| \cdot |G:M|.
\]
Exactly as in the prevous section, we see that $\Sigma_1 < q^{d-\frac{1}{2}}$ for large $q$.

Now let $M \in {\mathcal M}_2$. As in the previous proof we see that $|M|$ is larger than the size of the largest irreducible maximal subgroup of $G$ (which is at most $|Sp_n(q) \wr S_2|$ in the symplectic case, and at most $|GL_{n+1}(q).2|$ in the orthogonal case). Hence $M$ is in the class of reducible maximal subgroups of $G$. These are
\begin{itemize}
\item[(A)] stabilizers of non-degenerate subspaces, and
\item[(B)] parabolic subgroups.
\end{itemize}

Consider first case (A). For notational convenience, we deal with the case where $G = Sp_{2n}(q)$, the other cases being similar. In this case $M = Sp_{2r}(q) \times Sp_{2n-2r}(q)$ for some $r<\frac{n}{2}$. Since $n>N_5(\G) \ge 2N_4(\G)$ by hypothesis, we have $n-r > N_4(\G)$, and so Theorem \ref{fuchsclass}(ii) gives
\[
|{\rm Hom}(\G,M)| < g(n)\, q^d\,|Sp_{2n-2r}(q)|^{\mu+1}|Sp_{2r}(q)|^{vg+d-1}.
\]
Hence we see that  $|{\rm Hom}(\G,M)| \cdot |G:M| < q^E$, where
\[
E = \mu \dim M + \dim G + (\sum \frac{1}{m_i})\dim Sp_{2r}(q) + d.
\]
It follows that $E\le D-1$ if and only if the following inequality holds:
\[
\mu(\dim G - \dim M) > (\sum \frac{1}{m_i})(2r^2+r) + \frac{1}{2}\sum m_i^2 + d+1.
\]
Now $\dim G - \dim M = 4r(n-r)$ and $\mu > \sum \frac{1}{m_i}$ (by hypothesis), so the above inequality holds provided
\[
\mu(4rn-6r^2-r) >  \frac{1}{2}\sum m_i^2 + d+1.
\]
Since the left hand side is at least $\mu(4n-7)$, the inequality holds as
$$n>N_5(\G) > \frac{1}{2}\sum m_i^2 +2.$$
Hence $|{\rm Hom}(\G,M)| \cdot |G:M| < q^{D-1}$ for $M$ in case (A).

Now suppose $M$ is a maximal parabolic subgroup of $G$. Again, for convenience we just handle the case
$G = Sp_{2n}(q)$ and leave the very similar orthogonal cases to the reader. Let $M=P_r$, the stabilizer of a totally singular $r$-space, where $r\le n$. Then $M = QL$, with unipotent radical $Q$ and Levi subgroup
$L = GL_r(q) \times Sp_{2n-2r}(q)$. As in the previous section (\ref{rad}) holds.

Assume first that $r>N_2(\G)$ and $n-r>N_4(\G)$. Then Theorems \ref{fuchsian} and \ref{fuchsclass} imply that
$|{\rm Hom}(\G,L)| < g(n)q^{d+1}|L|^{\mu+1}$, and so by  (\ref{rad}),
\[
|{\rm Hom}(\G,M)| \cdot |G:M| < g(n)q^{d+1}|Q|^{vg+d}|L|^{\mu+1} < g(n) \,q^F,
\]
where $F = (vg+d)\dim Q+(\mu+1)\dim L + d+1$.
Using the fact that $\dim G = \dim L + 2\dim Q$, we find that $F\le D_1$ provided
\[
(\mu-\sum \frac{1}{m_i})\dim Q > \frac{1}{2}\sum m_i^2+d+2.
\]
Now $\dim Q = 2nr-\frac{3}{2}r^2+\frac{1}{2}r \ge \frac{1}{2}nr > \frac{1}{2}nN_2(\G)$. Since
$\mu-\sum \frac{1}{m_i}) \ge 1$ by hypothesis, the above inequality holds provided $nN_2(\G) > \sum m_i^2+2d+4$, which is true since $n>N_5(\G)$. Hence in this case,
$$|{\rm Hom}(\G,M)| \cdot |G:M| < q^{D-1}.$$

Now suppose that $r\le N_2(\G)$. Then as $n>N_5(\G) \ge 2N_2(\G),2N_4(\G)$, we have $n-r > N_4(\G)$, and so Theorem \ref{fuchsclass} gives
\[
|{\rm Hom}(\G,M)| \cdot |G:M| < g(n)q^{d}|Q|^{vg+d}|Sp_{2n-2r}(q)|^{\mu+1}|GL_r(q)|^{vg+d-1} < g(n) \,q^H,
\]
where $H = (vg+d)\dim Q + (\mu+1)\dim L +(\sum \frac{1}{m_i})r^2 + d$. As usual, we argue that $H<D-1$ provided
\begin{equation}\label{sill}
(\mu-\sum \frac{1}{m_i})\dim Q > (\sum \frac{1}{m_i})r^2+\frac{1}{2}\sum m_i^2+d+1.
\end{equation}
As $\dim Q = 2nr-\frac{3}{2}r^2+\frac{1}{2}r > nr$ (since $n>2N_2(\G)\ge 2r$), and also $\mu-\sum \frac{1}{m_i}) \ge 1$ by hypothesis, the above inequlaity holds provided
\[
nr > (\sum \frac{1}{m_i})r^2+\frac{1}{2}\sum m_i^2+d+1.
\]
This is true for $r\ge 2$ since $n>N_5(\G) \ge \s_1N_2(\G)+\s_3+2$ (notation as in the definition of $N_5(\G)$ before Theorem \ref{fuchsgenclass}); and for $r=1$ we check that (\ref{sill}) holds by putting in the exact value of $\dim Q$.
Hence in this case, we also have $|{\rm Hom}(\G,M)| \cdot |G:M| < q^{D-1}$.

Finally, assume that $n-r \le N_4(\G)$. Then $r>N_2(\G)$, so Theorem \ref{fuchsclass}(ii) gives
\[
|{\rm Hom}(\G,M)| \cdot |G:M| < g(n)q|Q|^{vg+d}|GL_r(q)|^{\mu+1}|Sp_{2n-2r}(q)|^{vg+d-1} < g(n) \,q^K,
\]
where $K = (vg+d)\dim Q + (\mu+1)\dim L +(\sum \frac{1}{m_i}\dim Sp_{2n-2r}(q) + 1$. As usual, $K<D-1$ provided
\[
(\mu-\sum \frac{1}{m_i})\dim Q > (\sum \frac{1}{m_i})\dim Sp_{2n-2r}(q)+\frac{1}{2}\sum m_i^2+2.
\]
Now $\dim Q \ge \frac{1}{2}(n^2+n)$ and $n-r\le N_4(\G)$, so the above inequality holds provided
\[
n^2+n > 2(\sum \frac{1}{m_i}) \dim Sp_{2N_4}(q) +  \sum m_i^2+4.
\]
This does hold, by the hypothesis $n > N_5(\G) \ge (1+2\s_1)N_4(\G)+\s_2$ (notation as in the definition of $N_5(\G)$.

We have now proved that $|{\rm Hom}(\G,M)| \cdot |G:M| < q^{D-1}$ for all the reducible maximal subgroups $M$. Hence the sum $\Sigma_2 < q^{D-\frac{1}{2}}$, and this completes the proof of Theorem \ref{fuchsgenclass}. \hal

\section{Probabilistic generation of exceptional groups}

In this section we prove Theorem \ref{excep2}.
 As in the previous sections, let $\G$ be a co-compact Fuchsian group of genus $g$ having $d$ elliptic generators $g_1,\ldots ,g_d$ of orders $m_1,\ldots ,m_d$, and define $\mu = \mu(\G)$ as in (\ref{meas}).

\subsection{A result for classical groups}\label{clpf}

For the proof of Theorem \ref{excep2}, we will need the following variant of Theorem \ref{excep1} for classical groups of small rank.

\begin{thm}\label{class1}
Let $\G$ be as above, and let $\GC$ be a simple adjoint algebraic group of type $A_r\,(r\le 7)$, $D_r\,(r\le 7)$,
$B_r\,(r\le 4)$ or $C_r\,(r\le 3)$ over an algebraically closed
field $\KK$ of good characteristic $p$. Suppose that $m_1\cdots m_d$ is coprime to $30$ and also to $p$.
 Then
\[
\dim {\rm Hom}(\G,\GC) = (vg-1)\dim \GC + \sum_{i=1}^d \dim J_{m_i}(\GC) +\d,
\]
where
\[
 0 \le \d \le  \left\{\begin{array}{l}
3, \hbox{ if }\GC = D_5, \,d=3 \hbox{ and }(m_1,m_2,m_3) = (7,7,7) \\
2, \hbox{ if }\GC = D_6 \hbox{ or }B_4, \,d=3 \hbox{ and }(m_1,m_2,m_3) = (7,7,7) \\
1, \hbox{ if }\GC = A_7 \hbox{ or }D_7, \,d=3 \hbox{ and }(m_1,m_2,m_3) = (7,7,7) \\
1, \hbox{ if }\GC = D_5, \,d=3 \hbox{ and }(m_1,m_2) = (7,7),\,m_3>7 \\
0, \hbox{ otherwise}.
\end{array}
\right.
\]
\end{thm}

The proof of this theorem runs along similar lines to that of Theorem \ref{excep1}, and rather than giving all the details we will sketch the proof for the case where $\GC = D_7$ and leave the other cases to the reader.

For the $D_7$ case we will need the bounds on $\a(\LC)$ for Levi subgroups provided by the next lemma.

\begin{lem}\label{alphad7}
Let $\GC=D_7$, and let $\LC$ be a Levi subgroup of $\GC$.
\begin{itemize}
\item[{\rm (i)}] Then $\a(\LC) < \frac{5}{6}$.
\item[{\rm (ii)}] If $\dim \LC \le 25$ then one of the following holds:
\begin{itemize}
\item[{\rm (a)}] $\LC \triangleright A_3$ and $\a(\LC) = \frac{1}{3}$;
\item[{\rm (b)}] $\a(\LC) < \frac{1}{4}$.
\end{itemize}
\end{itemize}
\end{lem}

\pf  (i) This follows from Theorem \ref{ratio}, apart from the case where $\LC' = D_6$. For this case it is routine to list the unipotent class representative $u$ of $\LC$ and check that $\frac{\dim u^\LC}{\dim u^\GC} < \frac{5}{6}$.

(ii) The Levi subgroups $\LC$ of $D_7$ of dimension 25 or less are those with $\LC' = A_1^k\,(k\le 4)$
$A_2A_1^k\,(k\le 3)$,  $A_3A_1^k\,(k\le 2)$, or $A_2^2$. In each case we list the unipotent class representatives of $\LC$ and check that conclusion (a) or (b) holds. \hal

\vspace{6mm}
\no {\bf Proof of Theorem \ref{class1} }
As stated above, we will just give the proof for the case where $\GC = D_7$. This follows along the lines of Section \ref{lines}.
Let $q$ be a power of $p$ such that $q \equiv 1 (\mod m_i)$ for all $i$, and let $F$ be a Frobenius endomorphism of $\GC$ such that $G = G(q) = \GC^F$ is a group of type $D_7$ over $\F_q$. For each $i$ define $J_{m_i}(\GC) = \{x \in \GC : x^{m_i}=1\}$, and let $j_{m_i}:=j_{m_i}(\GC) = \dim J_{m_i}(\GC)$.
We will prove that there are positive constants $c_1,c_2$ such that for sufficiently large $q$,
\begin{equation}\label{ulbds1}
c_1|G|^{vg-1}q^{\sum j_{m_i}} \le |{\rm Hom}(\G,G)| \le c_2|G|^{vg-1}q^{\sum j_{m_i}+\d},
\end{equation}
where $\d$ is as in Theorem \ref{class1}. Just as in previous proofs, this implies Theorem \ref{class1}.

Note that by \cite{law}, we have
\[
j_{m_i}(\GC) = 84-2\d_{m_i,11}-6\d_{m_i,7},
\]
and note also that $\frac{2}{h} = \frac{1}{6}$, where $h$ is the Coxeter number of $\GC$. For classes $x_i^G$ with $\dim x_i^\GC = j_{m_i}$, and $\CB_G(x_i) = L_i = \LC_i^F$, we have $\LC_i' = A_1^k$ with $k\le 3$, and it is easy to check thet $\a(\LC_i) \le \frac{1}{6}$ for all $i$. Hence we see as in Section \ref{lines} that (\ref{dec}) holds, giving the lower
bound in (\ref{ulbds1}).

Now we prove the upper bound.
For each $i$ let $y_i \in G$ be an element of order $m_i$, so that $\CB_G(y_i) = L_i = \LC_i^F$, where $\LC_i$ is an $F$-stable split Levi subgroup of $\GC$.
 Let $C_i = y_i^G$, ${\mathcal C}_i = y_i^\GC$
 and assume that $\dim {\mathcal C}_1 \le \cdots \le \dim {\mathcal C}_d$. Write ${\bf C} = (C_1,\ldots, C_d)$.
As in (\ref{zineq}) and (\ref{hombd}), we may assume that
\begin{equation}\label{zineq1}
d-2+vg-\sum_{i=1}^d \a(\LC_i) \le \frac{2}{h} = \frac{1}{6},
\end{equation}
and also, provided $d\ge 2$, we have
\begin{equation}\label{hombd1}
|{\rm Hom}_{\bf C}(\G,G)| \le c|G|^{vg-1}q^{\sum j_{m_i}}\cdot q^{\dim \GC + \dim {\mathcal C}_1-j_{m_1}-j_{m_d}}.
\end{equation}

If $d=0$ then (\ref{zineq1}) fails, and if $d=1$  then (\ref{zineq1}) gives $\a(\LC_1) \ge \frac{5}{6}$, contrary to Lemma \ref{alphad7}.

Now assume $d=2$, so that (\ref{zineq1}) gives $\a(\LC_1)+\a(\LC_2) \ge \frac{5}{6}$. Then Lemma \ref{alphad7} gives $\dim \LC_1 \ge 26$, so $\dim {\mathcal C}_1 \le 65$, which,  using (\ref{hombd1}), yields the upper bound in (\ref{ulbds1}).

Finally, suppose $d\ge 3$. By Lemma \ref{alphad7} we must have $\LC_i \triangleright A_3$ for $1\le i\le d-1$, since otherwise $\sum \a(\LC_i) < \frac{1}{3}(d-2)+\frac{1}{4}+\frac{1}{4}$, contrary to (\ref{zineq1}). Hence for $1\le i\le d-1$ we have $\dim {\mathcal C}_i \le \dim \GC - \dim A_3T_4 = 72$. It follows that
\begin{equation}\label{last1}
|{\rm Hom}_{\bf C}(\G,G)| \le c|G|^{vg}|C_1|\cdots |C_{d-1}| \le c|G|^{vg-1}q^{72(d-1)+91},
\end{equation}
which implies the upper bound in (\ref{ulbds1}). \hal

\subsection{Proof of Theorem \ref{excep2}}

This runs along similar lines to the proof of Theorem \ref{fuchsgensl} in Section \ref{slfuchs}.
Let $\GC$ be a simple adjoint algebraic group of exceptional type over an algebraically closed
field $K$ of good characteristic $p$. Suppose that $m_1\cdots m_d$ is coprime to $30$ and also to $p$. Let $q$ be a power of $p$ such that $q \equiv 1 \hbox { mod }m_i$ for all $i$, and let $F$ be a Frobenius endomorphism of $\GC$ such that $G = G(q) = \GC^F$ is a finite exceptional group of Lie type over $\F_q$. By (\ref{ulbds}) we have
$|{\rm Hom}(\G,G)| \sim q^D$, where
\[
D = (vg-1)\dim \GC + \sum_{i=1}^dj_{m_i}(\GC).
\]
Let $k$ be the lower bound for $j_{m_i}(\GC)$ in (\ref{jbds}), and define
\[
N = 3k-2\dim \GC -1,
\]
so that $N = 139$, 75, 41, 27 or 7, according as $G = E_8$, $E_7$, $E_6$, $F_4$ or $G_2$ respectively.

Now define ${\mathcal M}_1$ to be a set of conjugacy class representatives for the maximal subgroups $M$ of $G$ such that
$|M| < q^{N+\frac{1}{2}}$; define  ${\mathcal M}_2$ to be representatives for those $M$ such that $|M| \ge q^{N+\frac{1}{2}}$ and $M$ is non-parabolic; and define ${\mathcal M}_3$ to be a set of representatives of the maximal parabolic subgroups of $G$. For $i=1,2,3$ set
\[
\Sigma_i = \sum_{M \in {\mathcal M}_i} |{\rm Hom}(\G,M)| \cdot |G:M|.
\]
We aim to prove that there is a constant $c$ such that for large $q$,
\begin{equation}\label{comp}
\Sigma_1+ \Sigma_2 +\Sigma_3 < cq^{D-\frac{1}{2}}.
\end{equation}
This will complete the proof of Theorem \ref{excep2}.

First consider $\Sigma_1$. For $M \in {\mathcal M}_1$ we have
\[
|{\rm Hom}(\G,M)| \cdot |G:M| \le |M|^{vg+d-2}|G| \le cq^{N(vg+d-2)+\dim \GC}.
\]
From the definition of $N$ we check that $N(vg+d-2)+\dim \GC\le D-1$. Also $| {\mathcal M}_1| \le c'\log\log q$ by \cite[1.3]{LMS}, and it follows that for large $q$,
\begin{equation}\label{comp1}
\Sigma_1  < q^{D-\frac{1}{2}}.
\end{equation}

Now consider $\Sigma_2$. By \cite{LiSax}, ${\mathcal M}_2$ consists of the following subgroups:
\[
\begin{array}{|c|c|c|c|c|c|}
\hline
G & E_8(q) & E_7(q) & E_6^\e(q) & F_4(q) & G_2(q)  \\
\hline
{\mathcal M}_2 & - & \NB_G(^2\!E_6(q)) & \NB_G(D_5^\e(q)) & B_4(q) & \NB_G(A_2^\pm (q)) \\
\hline
\end{array}
\]
In all cases there are at most two conjugacy classes of the given subgroups. For $M \in {\mathcal M}_2$ we can use
(\ref{ulbds}) for $^2\!E_6(q)$, and (\ref{ulbds1}) for the other cases, to obtain
\[
|{\rm Hom}(\G,M)| \le cq^{(vg-1)\dim M+\sum j_{m_i}(M)+\d},
\]
where $\d$ is as in Theorem \ref{class1} (and as usual we write $\dim M$ for the degree of $|M|$ as a polynomial in $q$, and so on). Now compute that $(vg-2)\dim M + \sum j_{m_i}(M) + \dim G + \d \le D-1$. Hence
\begin{equation}\label{comp2}
\Sigma_2  < q^{D-\frac{1}{2}}.
\end{equation}

Finally, consider $\Sigma_3$. Let $M$ be a maximal parabolic subgroup of $G$, and wirte $M = QL$ where $Q$ is the unipotent radical and $L$ a Levi subgroup. As in (\ref{rad}) we have
\[
|{\rm Hom}(\G,QL)|\cdot |G:QL|  \le c|Q|^{vg+d}|{\rm Hom}(\G,L)|.
\]
Again using (\ref{ulbds}) and (\ref{ulbds1}),
\[
|{\rm Hom}(\G,L)| \le cq^{(vg-1)\dim L'+\sum j_{m_i}(L')+\d},
\]
where $\d \le 3$. Hence $|{\rm Hom}(\G,QL)| \cdot |G:QL| < cq^K$, where
\[
K = (vg+d)\dim Q+(vg-1)\dim L'+\sum j_{m_i}(L')+\d.
\]
Now using the fact that $\dim G = \dim L + 2\dim Q$, we check that $K\le D-1$, hence $\Sigma_3 < q^{D-\frac{1}{2}}$, if and only if the following inequality holds:
\begin{equation}\label{kine}
(d+2-vg)\dim Q \le \sum_{i=1}^d\left(j_{m_i}(G)-j_{m_i}(L')\right)+vg-2-\d.
\end{equation}
The values of $\dim Q$ can easily be computed for each maximal parabolic, and the values of $j_{m_i}(G)$ and $j_{m_i}(L')$ are given by \cite{law}. From this information we compute that the inequality (\ref{kine}) holds in all but the following case: $\GC = E_6$, $\LC' = D_5$ and $\d = 3$. Hence in all but possibly this exceptional case, we have shown that
\begin{equation}\label{comp3}
\Sigma_3  < q^{D-\frac{1}{2}}.
\end{equation}
Together with (\ref{comp1}) and (\ref{comp2}), this proves (\ref{comp}), completing the proof of Theorem \ref{excep2}.

To conclude the proof, we handle the exceptional case. The case $\d=3$ for $\LC' = D_5$ arises when $d=3$, $(m_1,m_2,m_3) = (7,7,7)$, in which case the argument for (\ref{last1}) only gives
\[
|{\rm Hom}_{\bf C}(\G,L)| \le c|L'|^{vg}|C_1|\,|C_2| \le c|L'|^{vg-1}q^{36\cdot 2+45},
\]
where ${\bf C} = (C_1,C_2,C_3)$ and  $C_i = y_i^{L'}$ ($i=1,2,3$) are classes in $L'$ for which $\CB_{\LC}(y_i) = A_1^2T_3$. In the parabolic $QL$ of $G = E_6^\e(q)$, the unipotent radical $Q$ is a spin module under the action of $L$, and we compute that $\dim \CB_Q(y_i) = 4$. Hence
\[
|{\rm Hom}_{\bf C}(\G,QL)|\cdot |G:QL| \le c|QL|^{vg}|y_1^{QL}| \,|y_2^{QL}|\cdot |Q| \le c'|QL|^{vg}q^{114},
\]
and so the contribution $|{\rm Hom}_{\bf C}(\G,QL)|\cdot |G:QL|$ of the class triple ${\bf C}$ to the sum $\Sigma_3$ is less than $|G|^{vg}q^{\sum j_{m_i}(G)-\dim G} = |G|^{vg}q^{66\cdot 3-78}$. Hence again
$\Sigma_3  < q^{D-\frac{1}{2}}$, completing the proof of the theorem. \hal


\end{document}